\documentclass[reqno]{amsart}
\usepackage{amssymb,amscd}
\newcommand{\Wcoho}{H_{\rm QDS}}

\newcommand{\Fdash}{{F'}}

\newcommand{\B}{{\mathcal{B}}}

\newcommand{\N}{{\mathcal{N}}}

\newcommand{\DW}{{\bf D}^{\mathcal W}}

\newcommand{\Dtot}{\widehat{\bf D}}
\newcommand{\Dg}{{\bf D}}

\newcommand{\isomap}{\overset{\sim}{\rightarrow} }
\newcommand{\W}{{\mathcal{W}}}
\newcommand{\Whw}[2]{{\bf{h}}_{#1}^{#2}}
\renewcommand{\t}{{\mathfrak{t}}}

\newcommand{\sI}{\bar{I}}
\newcommand{\1}{{\mathbf{1}}}
\newcommand{\Ln}[1]{L\bar{\mathfrak{n}}_{#1}}

\newcommand{\fd}{{\mathcal{C}}}

\newcommand{\w}{{\textsl{w}}}
\newcommand{\sQ}{\bar{Q}}
\newcommand{\eW}{\widetilde{W}}
\renewcommand{\Pr}{\textsl{{Pr}}_{\kappa}}

\newcommand{\sroots}{\bar{\Delta}}
\newcommand{\roots}{\Delta}
\newcommand{\proots}{\Delta_+}

\newcommand{\rroots}{\Delta^{\rm re}}
\newcommand{\iroots}{\Delta^{\rm im}}
\newcommand{\prroots}{\Delta_+^{\rm re}}
\newcommand{\piroots}{\Delta_+^{\rm im}}
\newcommand{\nrroots}{\Delta_-^{\rm re}}
\newcommand{\F}{{\mathcal{F}}}

\newcommand{\m}{{\mathfrak{m}}}

\renewcommand{\a}{{\mathfrak{a}}}

\renewcommand{\a}{{\mathfrak{a}}}

\newcommand{\sPi}{\bar{\Pi}}

\newcommand{\Cl}{{\mathcal{C}l}}

\newcommand{\sproots}{\bar{\Delta}_+}
\newcommand{\snroots}{\bar{\Delta}_-}

\newcommand{\coho}[1]{H^{\frac{\infty}{2}+{#1}}}

\newcommand{\Hsim}[1]{I({#1})}
\newcommand{\semiinf}{\frac{\infty}{2}}
\renewcommand{\H}{{\widetilde{\h}}}

\newcommand{\C}{{\mathbb C}}
\newcommand{\Z}{{\mathbb Z}}
\newcommand{\Q}{{\mathbb Q}}
\newcommand{\inv}{^{-1}}
\newcommand{\dual}[1]{{#1}^*}
\newcommand{\lam}{\lambda}
\newcommand{\Lam}{\Lambda}

\renewcommand{\*}{{\otimes}}
\newcommand{\+}{\mathop{\oplus}}
\newcommand{\h}{ {\mathfrak h}}
\newcommand{\g}{ {\mathfrak g}}
\renewcommand{\b}{{\mathfrak{b}}}

\newcommand{\che}{^{\vee}}

\newcommand{\bra}{{\langle}}
\newcommand{\ket}{{\rangle}}

\DeclareMathOperator{\tr}{tr}

\DeclareMathOperator{\End}{End}

\DeclareMathOperator{\Hom}{Hom}

\DeclareMathOperator{\id}{id}
\DeclareMathOperator{\ad}{ad}
\DeclareMathOperator{\haru}{span}

\DeclareMathOperator{\Ext}{Ext}

\DeclareMathOperator{\ch}{ch}

\DeclareMathOperator{\rank}{rank}

\DeclareMathOperator{\height}{ht}

\newcommand{\sg}{ \bar {\mathfrak g}}
\newcommand{\sh}{\bar \h}
\newcommand{\sn}{\bar{\mathfrak{n}}}
\newcommand{\snp}{\bar{\mathfrak{n}}_+}
\newcommand{\snn}{\bar{\mathfrak{n}}_-}

\newcommand{\sP}{\bar P}
\newcommand{\sW}{\bar W}
\newcommand{\srho}{\bar \rho}

\renewcommand{\srho}{\bar \rho}

\renewcommand{\check}{^{\vee}}

\DeclareMathOperator{\im}{Im}

\newcommand{\ud}[2]{{\genfrac{}{}{0pt}{}{#1}{#2}}}

\newcommand{\n}{{{\mathfrak{n}}}}

\newcommand{\BGG}{{\mathcal O}}

\numberwithin{equation}{section}

\title[Vanishing of
cohomology
associated to
quantized  reduction]{Vanishing of
cohomology
associated to
\\quantized Drinfeld-Sokolov reduction}
\author{Tomoyuki Arakawa}
\address{Graduate school of Mathematics,  Nagoya
University,
Chikusa-ku, Nagoya, 464-8602, JAPAN}
\email{tarakawa@math.nagoya-u.ac.jp}
\begin{document}
%\today
%%%%%%%% Thoremstyles %%%%%%%%%%
\theoremstyle{plain}
\newtheorem{Th}{Theorem}[section]
\newtheorem{MainTh}{Main theorem}
\newtheorem{Pro}[Th]{Proposition}
\newtheorem{Lem}[Th]{Lemma}
\newtheorem{Co}[Th]{Corollary}
\renewcommand{\theMainTh}{}

\newtheorem{Facts}[Th]{Facts}

\theoremstyle{definition}

\newtheorem{dfandpr}[Th]{Definition and Proposition}
\theoremstyle{remark}
\newtheorem{Def}[Th]{Definition}
\newtheorem{Rem}[Th]{Remark}
\newtheorem{Conj}{Conjecture}
\newtheorem{Claim}{Claim}
\newtheorem{Notation}{Notation}
\newtheorem{Ex}[Th]{Example}

\newcommand{\st}{{\mathrm{st}}}
%\numberwithin{equation}{subsection}
%personal preprint, \today
\maketitle
\begin{abstract}
We prove a  vanishing
theorem
of the cohomology
arising from the
two quantized Drinfeld-Sokolov reductions
(^^ ^^ $+$" and ^^ ^^ $-$" reduction)
introduced by
Feigin-Frenkel and Frenkel-Kac-Wakimoto.
As a consequence,
the  vanishing conjecture
of Frenkel-Kac-Wakimoto
is proved
for the ^^ ^^ $-$" reduction
and partially for  the ^^ ^^ $+$" reduction.
\end{abstract}
\section{Introduction}
%\subsection{}
In this paper
we study the
cohomology
of the BRST complex of the
quantized Drinfeld-Sokolov reductions
%at non-critical level
introduced by Feigin-Frenkel \cite{FF_W}
and Frenkel-Kac-Wakimoto \cite{FKW}
in their study of $W$-algebras.

\medskip

Let $\sg=\snn\+ \sh \+\snp$ be a
finite-dimensional complex simple Lie algebra.
Let
$\g=\sg\* \C[t,t\inv]
\+\C K\+ \C \Dg$ be the affine Lie algebra
associated to
$\sg$.
Let
$\prroots$ be the set of real positive roots
of $\g$,
$\sproots\subset \prroots$ the set of positive roots
of $\sg$,
$W$ the Weyl group of $\g$.
Let $\kappa\in \C\backslash \{0\}$
and let $\dual{\h}_{\kappa}$
denote
the
set of the
weights of $\g$ of level
$\kappa-h\che$.
Let $\BGG_{\kappa}$ be
the Bernstein-Gelfand-Gelfand category
of $\g$ of level $\kappa-h\che$,
where
$h\check$ is the dual Coxeter number of $\sg$.
Let $L(\Lam)$,
$\Lam\in \dual{\h}_{\kappa}$,
be the simple module of $\BGG_{\kappa}$
of highest weight $\Lam$.

Let $\Ln{\pm}=\sn_{\pm}\*
\C[t,t\inv]\subset
\g$.
Fix  a nondegenerate character
$\bar{\chi}_{\pm}$  (\cite{Kostant-Whittaker})
of
$\sn_{\pm}$ as in \cite[2.1]{FKW}.
It extends to a character
${\chi}_{\pm}:
\Ln{\pm}\rightarrow \C$
by
\begin{align*}
&\chi_{+}(X\*t^n)=\delta_{n,-1}\bar{\chi}_+(X)
\quad (X\in \snp,n\in \Z),\\
&\chi_{-}(X\*t^n)=\delta_{n,0}\bar{\chi}_-(X)
\quad (X\in \snn,n\in \Z).
\end{align*}
Let
$\C_{\chi_{\pm}}$ be the
one-dimensional
representation
of $U(\Ln{\pm})$
defined by $\chi_{\pm}$.
Then,
the semi-infinite cohomology
$\Wcoho^{\bullet}(\Ln{\pm},V)
=\coho{\bullet}(\Ln{\pm},V\* \C_{\chi_{\pm}})$,
$V\in \BGG_{\kappa}$,
is called
cohomology
of the BRST complex of the
quantized Drinfeld-Sokolov reduction
for $\Ln{\pm}$
(^^ ^^ $+$" and ^^ ^^ $-$" reduction)
associated to  $V$
(\cite{FF_W,FKW,FB}).

Let $V_{\kappa}(\sg)
=U(\g)\*_{U(\sg\* \C[t]\+\C K
\+\C \Dg) }\C$
be the
universal affine vertex algebra
associated to $\sg$
of
level
$\kappa-h\che$.
Then,
the $0$-th cohomology
$\Wcoho^{0}(\Ln{+},V_{\kappa}(\sg))$
is the Feigin-Frenkel's $W$-algebra
$\W_{\kappa}(\sg)$ associated to $\g$
of level $\kappa-h\che$
(\cite{FF_W}).
Their realization of
$\W_{\kappa}(\sg)$
gives a  functor
\begin{align}\label{eq:FF-functor}
V \rightsquigarrow
\Wcoho^{i}(\Ln{\pm},V)\quad
(i\in \Z)
\end{align}
from
$\BGG_{\kappa}$
to
the category of
$\W_{\kappa}(\sg)$-modules
(\cite{FF_W,FKW,FB}).

\medskip

Let us now describe our result.
Let
$W^{\Lam}\subset W$ be the
integral Weyl group
of a weight $\Lam\in \dual{\h}$.
For  $\Lam\in \dual{\h}_{\kappa}$,
let $\BGG_{\kappa}^{[\Lam]}$
be the full subcategory
of $\BGG_{\kappa}$
whose objects have all their
local composition factors
isomorphic to $L(w\circ\Lam)$,
%the simple module of highest
%weight
%$w(\Lam+\rho)-\rho$,
$w\in W^{\Lam}$.
Then,
$
\BGG_{\kappa}=\bigoplus\limits_{\Lam
\in \dual{\h}_{\kappa}/\sim}\BGG_{\kappa}^{[\Lam]}
$,
where
$\sim$ is the equivalent relation
defined by
$\lam\sim \mu\Leftrightarrow
\mu \in W^{\lam}\circ \lam$.
The main result of this paper is the following.
\begin{MainTh}
Let $\Lam\in \dual{\h}_{\kappa}$,
$\kappa\in \C\backslash\{0\}$.
Then,
\begin{enumerate}
\item $\Wcoho^{i}(\Ln{-},V)=\{0\}$
$(i\ne 0)$
for all objects $V$
in $\BGG^{[\Lam]}_{\kappa}$
if $\bra \Lam,\bar{\alpha}\che\ket\not\in \Z
$ for all $\bar{\alpha}\in \sproots$,
\item
$\Wcoho^{i}(\Ln{+},V)=\{0\}$
$(i\ne 0)$
for all objects $V$
in $\BGG^{[\Lam]}_{\kappa}$
if $\bra \Lam,\alpha\che\ket\not\in \Z$
for all
$\alpha\in \prroots\cap t_{\srho\che}(\nrroots)
=\{-\bar{\alpha}+ n\delta;
\bar{\alpha}\in \sproots,1\leq n\leq \height \bar{\alpha}\}$.
\end{enumerate}
\end{MainTh}
This result appears as Theorem
\ref{Theorem:Main-Theorem} in this paper.
This
shows that
the correspondence
$V\rightsquigarrow \Wcoho^{0}(L\sn_{\pm},V)$
defines an exact functor
from $\BGG_{\kappa}^{[\Lam]}$
to the category of $\W_{\kappa}(\sg)$-modules
under the condition of $\Lam$ described  as
above.
%This gives the character of
%$\Wcoho^0(L\sn_{\pm},L(\lam))$.
The irreducibility
of
$\Wcoho^0(L\sn_{\pm},L(\lam))$
will be studied
in our forthcoming
paper.

\medskip

Frenkel-Kac-Wakimoto
\cite{FKW} applied the functor
\eqref{eq:FF-functor}
to the  principal admissible
representations of $\g$
of fractional levels $\kappa-h\che$.
They conjectured
that
a
vanishing of cohomology
holds for that case
and that
the functor
\eqref{eq:FF-functor} sends a
principal admissible
representation to zero
or to an irreducible ^^ ^^ minimal"
representations
of $\W_{\kappa}(\sg)$.
Based on the vanishing
conjecture,
they
calculated the
characters and
fusion coefficients for
conjectural ^^ ^^ minimal"
representations of $\W_{\kappa}(\g)$.
Our result
settles the
vanishing conjecture
of Frenkel-Kac-Wakimoto
for
the ^^ ^^ $-$" reduction
and partially for  the ^^ ^^ $+$" reduction.
Though our result
for the ^^ ^^ $+$" reduction
is partial,
we  remark that
every conjectural irreducible
^^ ^^ minimal''
representation
is isomorphic to
$\Wcoho^{0}(\Ln{+},L(\Lam))
$
for some
principal admissible weight
$\Lam$
which satisfies the
condition of Main theorem
(2),
see Remark \ref{Rem:final} \eqref{rem-on-minimal}.

\medskip

%\subsection{}
This article is organized as follows.
In section \ref{section:Preliminaries},
we collect the necessary information
about the
affine Lie algebra $\g$
and its representations.
In section \ref{section:freeness},
we prove Theorem \ref{Th:cofreeness_of_Verma}
and Theorem \ref{Th:global_cofreeness}
which
is needed in the later arguments.
In section \ref{section:complex},
we recall the definition
of the cohomology $\Wcoho^{\bullet}(\Ln{\pm},V)$
and define some operators acting on the
corresponding complex.
In particular,
we define the
degree
operator
which acts on $\Wcoho^{\bullet}(\Ln{\pm},V)$,
$V\in \BGG_{\kappa}$,
semisimply.
This is essentially
the operator
$-L_0^{\pm}$ defined in
\cite[3.1]{FKW}.
The difficulty
dealing with this cohomology
arises from the fact that
by construction the corresponding
eigenspaces of complexes themselves are not
finite-dimensional in general.
The results in section
\ref{section:vanishing_of_projectives}
are straightforward generalization
of
\cite[14.2]{FB}.
Thus,
Theorem \ref{Th:vanishing.for.Verma},
which states the vanishing
of the cohomology
with coefficient
in Verma modules,
was essentially proved in \cite{FB}.
In  section \ref{section:vanishing_of_dual_of_Verma},
we prove
the corresponding statement for the
dual $M(\lam)^*$ of Verma
module $M(\lam)$ (Theorem \ref{Th:vanishin-dual}).
Though the usual duality (\cite{Feigin})
of semi-infinite cohomology
cannot be applied
for this cohomology,
this is done
by establishing the duality
\begin{align*}
&\Wcoho^{i}(\Ln{+},M(\lam)^*)
\cong
\Wcoho^{-i}(\Ln{-},M(t_{-\srho\che}\circ\lam))^*,\\
&\Wcoho^{i}(\Ln{-},M(\lam)^*)
\cong
\Wcoho^{-i}(\Ln{-},M(w_0\circ\lam))^*,
\end{align*}
for $i\in \Z$
under the similar
restriction of $\lam$
as in the main theorem.
Here,
${}^*$ is the graded dual
and $w_0$ is the longest element of the Weyl group
of $\sg$.
This result  is proved by
using some spectral sequences.
The duality above may
explain why the ^^ ^^ $-$" reduction
behaves ^^ ^^ nicer".
In section \ref{section:estimate},
we estimate the eigenvalues
of degree operators on $\Wcoho^{\bullet}(\Ln{\pm},V)$.
%(Proposition \ref{Pro:pre-crutial-estimate}).
The results in this section
play a crucial rule
in proving our main theorem
when $\kappa\in
\Q_{>0}$,
that is,
when the objects in $\BGG_{\kappa}$
do not necessarily have finite length.
Finally,
we give
the proof of
our main theorem
in
section \ref{section:Main-theorem}.

\bigskip

{\em Acknowledgments. }
I wish to thank Akihiro Tsuchiya
for his continuous encouragement,
Kiyokazu Nagatomo for his interest
and encouragement,
Kenji Iohara for useful discussion,
Hiroyuki Ochiai for
reading the manuscript.
I am grateful to Edward Frenkel for his
interests and valuable discussions.

\bigskip
\section{
Preliminaries }
\label{section:Preliminaries}
\subsection{Affine Lie algebra}
In the sequel,
we fix a nonzero complex number $\kappa$, %$\kappa\in \C\backslash \{0\}$,
a simple
finite-dimensional complex Lie algebra $\sg$ and a
Cartan subalgebra $\sh$. Let $\sroots$ denote the set
of roots,
$\sPi$  a basis of $\sroots$,
$\sproots$  the set of positive roots,
and $\snroots=-\sproots$.
This gives the triangular
decomposition $\sg=\snn\+\sh\+\snp$.
Let $\sQ$ denote the root lattice,
$\sP$ the weight lattice,
$\sQ\che$ the coroot lattice
and
$\sP\che$ the coweight lattice.
%$\sQ_+=\sum_{\alpha\in \sproots}\Z_{\geq 0}
%\alpha\subset \sQ$.
Let $\srho$ be the half  sum of
positive roots,
$\srho\che$ the half  sum of positive coroots.
For $\alpha\in \sproots$,
the number
$\bra \alpha,\srho\che\ket$
is called the
{\em height} of
$\alpha$ and denote by
$\height \alpha$.
Let $\sW$ be the Weyl group of $\sg$,
$w_0$ the longest element of $\sW$.

Let $(~,~)$ be the normalized invariant inner product
of $\sg$.
Thus,
$(~,~)=\frac{1}{2h\check}$Killing form,
where
$h\check$ is the dual Coxeter number of $\sg$.
We identify $\sh$ and $\dual{\sh}$ using
the form.
Then,
$\alpha\che
=2\alpha/(\alpha,\alpha)$,
$\alpha\in \sroots$.

\smallskip

Let
$\g=\sg\*\C[t,t^{-1}]\+\C K \+ \C \Dg$
be the affine Lie algebra associated to
($\sg$,$(~,~)$),
where
$K$ is its central element and $\Dg$ is the degree
operator
(see \cite{KacBook}).
The bilinear form $(~,~)$
is naturally extended from $\sg$ to $\g$.
Set $X(n)=X\* t^n$,
$X\in \sg$,
$n\in \Z$.
The subalgebra $\sg\* \C\subset \g$
is naturally identified with $\sg$.

Fix
the  triangular decomposition
$\g=\g_-\+\h\+\g_+$
in the standard way.
Thus,
\begin{align*}
&\h=\sh\+\C K\+\C \Dg,\\
&\text{$\g_-=\snn\* \C[t\inv]\+ \sh\* \C[t\inv]t\inv
\+ \snp\*\C[t\inv]t\inv$,}\\
&\text{$\g_+=\snn\* \C[t]t\+ \sh\* \C[t]t
\+ \snp\*\C[t]$.}
\end{align*}
Let
$\dual{\h}=\dual{\sh}\+\C \Lam_0\+\C \delta$
be the  dual of $\h$.
Here,
$\Lam_0$ and
$\delta$ are dual elements of $K$ and $\Dg$
respectively. For $\lam\in\dual{\h}$,
the number
$\bra \lam,K\ket$ is called the
{\em level of $\lam$}.
Let
$\dual{\h}_{\kappa}$
denote the
set of the
weights of level
$\kappa-h\che$:
\begin{align*}
\dual{\h}_{\kappa}
=\{\lam\in \dual{\h};\bra \lam+\rho,K\ket=\kappa\},
\end{align*}
where,
$\rho=\bar{\rho}+h\che
\Lam_0\in\dual{\h}$.
%$$\dual{\h}_{\kappa}=\{\lam\in\dual{\h}
%\mid \bra \lam+\rho,K\ket=\kappa\}.$$
Let
$\bar{\lam}$ be the restriction of $\lam\in \dual{\h}$ to $\dual{\sh}$.

\smallskip

Let $\roots$ be the set of roots of $\g$,
$\roots_+$ the set of positive roots,
$\roots_-=-\roots_+$.
Then,
$\roots=\rroots\sqcup \iroots$,
where
$\rroots$
is
the set of real roots
and
$\iroots$ is the set of imaginary roots.
Let
$\Pi$ be the basis of $\rroots$,
$\rroots_{\pm}=\rroots\cap \roots_{\pm}$,
$\iroots_{\pm}=\iroots\cap \roots_{\pm}$.
Let $Q$ be the root lattice,
$Q_+=\sum_{\alpha\in \proots}\Z_{\geq 0} \alpha\subset
Q$.

Let
$W\subset GL(\dual{\h})$ be the Weyl group of $\g$
generated by the reflections $s_{\alpha}$,
$\alpha\in \rroots$,
defined by
$s_{\alpha}(\lam)=\lam-\bra \lam,\alpha\che\ket\alpha$.
Then,
$W=\sW\ltimes \sQ\che$.
Let $\eW=\sW\ltimes \sP\che$,
the extended Weyl group of $\g$.
For $\mu\in \sP\che$,
we denote the corresponding element
of $\eW$ by $t_{\mu}$.
Then,
\begin{align*}
t_{\mu}(\lam)=\lam+\bra\lam,K\ket\mu-
\left(\bra\lam,\mu\ket+\frac{1}{2}
|\mu|^2\bra\lam,K\ket
\right)\delta \quad (\lam\in \dual{\h}).
\end{align*}
Let $\eW_+=\{w\in \eW;\prroots\cap
w\inv(\nrroots)=\emptyset  \}$.
Then,
$\eW=\eW_+\ltimes W $.

The dot action of
$\eW$
on $\dual{\h}$ is defined by $w\circ \lam= w(\lam+\rho)-\rho$
($\lam\in \dual{\h}$).

\smallskip

For $\Lam\in \dual{\h}$,
let
$R^{\Lam}=\{
\alpha\in \rroots;
\bra \Lam+\rho,\alpha\che\ket\in \Z\}$,
$R^{\Lam}_+=R^{\Lam}\cap
\prroots$,
$\Pi^{\Lam}=\{\alpha\in R^{\Lam}_+;
s_{\alpha}(R^{\Lam}_+\backslash \{\alpha\})
\subset R^{\Lam}_+\}$.
It is known that
$R^{\Lam}$ is a subroot
system of $\rroots$
with the basis $\Pi^{\Lam}$
(\cite{Moody,KT1}).
Let $Q^{\Lam}=\sum_{\alpha\in
\Pi^{\Lam}}\Z \alpha \subset Q$,
$Q^{\Lam}_+=\sum_{\alpha\in
\Pi^{\Lam}}\Z_{\geq 0} \alpha$.
For $\mu=\sum_{\alpha\in \Pi^{\Lam}}
m_{\alpha}\alpha
\in Q^{\Lam}_+$,
$m_{\alpha}\in \Z_{\geq 0}$,
set
$\height_{\Lam} (\mu)=\sum_{\alpha\in \Pi^{\Lam}}
m_{\alpha}$. Let $W^{\Lam}
=\bra s_{\alpha};\alpha\in R^{\Lam}\ket$
be the integral Weyl group corresponding
to $\Lam$.
Then,
$
R^{w\circ \Lam}=R^{\Lam}$
for $w\in W^{\Lam}$.

\subsection{BGG category
of $\g$}
For a $\g$-module $V$
(or for  simply a $\h$-module $V$),
let $V^{\lam}=\{v\in V;
hv=\lam(h)v \text{ for }h\in \h\}$
be the weight space of weight $\lam$.
Let $P(V)=\{\lam\in \dual{\h};V^{\lam}\not=
\{0\}\}$.
If
$\dim V^{\lam}<\infty$ for all $\lam$,
then
we set
\begin{align}\label{eq:graded-dual}
V^*=\bigoplus_{\lam}\Hom_{\C}(V^{\lam},\C)
\subset \Hom_{\C}(V,\C).
\end{align}
The formal character
$\ch V$ of
$V$ is defined as $\ch V=
\sum_{\lam}e^{\lam}\dim_{\C}V^{\lam}$.

\smallskip

Let
$\BGG_{\kappa}$ be the
full subcategory of the category of left $\g$-modules
consisting of objects $V$
such that
(1) $V$ is locally finite
over $\g_+$,
(2) $V=\bigoplus\limits_{\lam\in
\dual{\h}_{\kappa}}V^{\lam}$ and
$\dim_{\C}V^{\lam}<\infty$ for all $\lam$, (3)
there exists
a finite subset
$\{\mu_1,\dots,\mu_n\}\subset
\dual{\h}_{\kappa}$
such that
$P(V)\subset  \bigcup\limits_i \mu_i-Q_+$.

The correspondence
$V\rightsquigarrow V^*$
defines
the duality functor
in $\BGG_{\kappa}$.
Here, $\g$ acts on $V^*$
by
$(Xf)(v)=f(X^tv)$,
where
$X\mapsto X^t$ is  the Chevalley
antiautomorphism.
For a subalgebra $\a\subset \g$,
we set
\begin{align*}
\a^t=\{X^t;X\in \a\}\subset \g.
\end{align*}

Let $M(\lam)\in \BGG_{\kappa}$,
$\lam\in \dual{\h}_{\kappa}$,
be the Verma module of
highest weight $\lam$
and
$L(\lam)$ its unique simple
quotient.
Let $\BGG_{\kappa}^{[\Lam]}$,
$\Lam\in \dual{\h}_{\kappa}$,
be the
full subcategory of
$\BGG_{\kappa}$
whose objects have all their
local composition factors
isomorphic to $L(w\circ\Lam)$,
$w\in W^{\Lam}$.
By \cite{Kumer},
$\BGG_{\kappa}$ splits into the orthogonal
direct sum
$
\BGG_{\kappa}=\bigoplus\limits_{\Lam
\in \dual{\h}_{\kappa}/\sim}\BGG_{\kappa}^{[\Lam]}
$,
where
$\sim$ is the equivalent relation
defined by
$\lam\sim \mu\Leftrightarrow
\mu \in W^{\lam}\circ \lam$.
Orthogonal here means that
$\Ext^i_{\BGG_{\kappa}}(M,N)=0 $
for $M\in \BGG_{\kappa}^{[\Lam]}$,
$N\in \BGG_{\kappa}^{[\Lam']}$,
$i\geq 0$,
when $\Lam\ne \Lam'$ in $\dual{\h}_{\kappa}/\sim$.
%Then,
%$V^*\in \BGG_{\kappa}^{[\leq \Lam]}$
%for all $V\in \BGG_{\kappa}^{[\leq \Lam]}$.
\section{Some results on  $\BGG_{\kappa}$}
\label{section:freeness}
\subsection{}
For $w\in \eW$,
let $\g_w=\g_+\cap w(\g_-)\subset \g_+$.
Then,
$\g_w^t=\g_-\cap w(\g_+)\subset \g_-$.
In this section
we shall prove the following two
theorems which will be needed in the later arguments.

\begin{Th}\label{Th:cofreeness_of_Verma}
Let $\lam\in \dual{\h}$.
Suppose that
$\bra \lam+\rho,\alpha\che\ket\not\in \Z_{\geq 1}$
for all
$\alpha\in \prroots\cap w(\nrroots)$.
Then,
$M(\lam)$ is cofree over
$\g_w$.
\end{Th}
\begin{Th}\label{Th:global_cofreeness}
Let $w\in \eW$
and
$\Lam\in \dual{\h}_{\kappa}$
such that
$\bra \Lam+\rho,\alpha\che\ket\not\in \Z$
for all $\alpha\in
\prroots\cap w
(\nrroots) $.
Then,
%any object in $\BGG_{\kappa}^{[\Lam]}$
%is cofree over $\g_w$,
%or equivalently,
any object in $\BGG_{\kappa}^{[\Lam]}$
is free over $\g_w^t$.
\end{Th}

\subsection{}
Let us start with the following lemma:
\begin{Lem}\label{Lem:supplimentary}
Let $\m$ be an $\ad\h$-stable subalgebra of $\g_-$.
Let $V$ be a module over $\m\+\h\subset \g$
such that
$V=\bigoplus_{\lam\in \dual{\h}}V^{\lam}$
and $P(V)\subset  \bigcup\limits_i \mu_i-Q_+$
for some
finite subset
$\{\mu_1,\dots,\mu_n\}\subset
\dual{\h}$.
Suppose that
$V=\m V$.
Then,
$V=\{0\}$.
\end{Lem}
\begin{proof}
Suppose that
$V\ne \{0\}$.
Then,
there exists $\mu\in P(V)$
such that
$\mu+ \alpha\not\in P(V)$
for any $\alpha\in Q_+$.
But this contradicts $V=\m V$.
\end{proof}
\begin{Pro}\label{Pro:criterion-of-freeness}
Let $\m$ be an $\ad\h$-stable subalgebra of $\g_-$.
Then,
for $V\in \BGG_{\kappa}$,
the following conditions are equivalent:
\begin{enumerate}
\item $V$ is free over $\m$.
\item $H_1(\m,V)=0$.
\end{enumerate}
\end{Pro}
\begin{proof}
Clearly (1) implies (2).
Let us show (2) $\Rightarrow$ (1).
Let $\{\bar{v_j}; j\in J\}$ be a
basis of $H_0(\m,V)=V/\m V$
and let $v_j$, $j\in J$, be  an inverse image of
$\bar{v_j}$ in $V$.
Since $V/\m V$
is naturally a $\h$-module,
we may suppose that
each $v_j$ is a weight vector of $V$.
We claim that
$%\begin{align}%\label{eq:suppl1}
V=\sum_{j\in J}U(\m)v_j
$.
%\end{align}
Indeed,
$V=\sum_{j\in J}U(\m)v_j+ \m V$.
Let $\bar{V}=V/\sum_{j\in J}U(\m)v_j$.
Then,
$\bar{V}=\m\bar{V}$
and $P(\bar{V})\subset P(V)$.
Thus,
$\bar{V}=\{0\}$
by Lemma \ref{Lem:supplimentary}.

Let $V_1=\bigoplus_{j\in
J}U(\m)v_j$,
the free $U(\m)$-module with
a basis
$\{v_j; j\in J\}$.
Let $\h$ act semisimply on
$V_1$ so that
the natural map
$\pi:V_1\twoheadrightarrow V$
is a homomorphism of $(\m\+\h)$-module.
Let $M=\ker \pi$.
%Since
%$V_1=\bigoplus_{\lam} V_1^{\lam}$
%and $P(V_1)\subset P(V)+P(U(\m))$,
%it follows that
Then,
$M=\bigoplus_{\lam} M^{\lam}$
and $P(M)\subset P(V)+P(U(\m))$.
Here,
$\h$ acts on $U(\m)$
by adjoint.
Now suppose that $H_1(\m,V)=0$.
Then,
by the long exact sequence
\begin{align*}
\dots \rightarrow
H_1(\m,V)\rightarrow H_0(\m,M)\rightarrow
H_0(\m,V_1)\rightarrow H_0(\m,V)\rightarrow 0
\end{align*}of $\m$-homology,
it follows that
$H_0(\m,M)=0$,
that is,
$M=\m M$.
Hence
$M=\{0\}$
by
Lemma \ref{Lem:supplimentary}.
%Then,
%$M=\bigoplus_{\lam\in \dual{\h}}M^{\lam}$
%and $P(M)\subset  \bigcup\limits_i \mu_i-Q_+$
%for some
% finite subset
%$\{\mu_1,\dots,\mu_n\}\subset
%\dual{\h}_{\kappa}$.
%But
%Hence
% $M=\{0\}$
%by Lemma \ref{Lem:supplimentary}.
\end{proof}
\subsection{Arkhipov's  twisting functor
}
%Let
%$U(\g_w)^*=\bigoplus_{\lam}\Hom(U(\g_w)^{\lam},\C)$,
%where
%$U(\g_w)^{\lam}$ is the weight
%space of  $U(\g_w)$ with respect to
% the adjoint acion of $\h$.
Let $S_w$
be the
Arkhipov's
semiregular module
corresponding to $w\in W$ (\cite{Ark},
see also \cite{AL,S}).
It is a $U(\g)$-bimodule
and
\begin{align*}
S_w&=U(\g_w)^*\*_{U(\g_w)}U(\g)
\quad \text{(as a left $U(\g_w)$-module and a right
$U(\g)$-module) }\\
&=U(\g)\*_{U(\g_w)}U(\g_w)^*
\quad \text{(as a left $U(\g)$-module
and a right $U(\g_w)$-module)}.
\end{align*}
Here,
$U(\g_w)^*$ is considered to
be  a $U(\g_w)$-bimodule
by $(fx)(n)=f(xn)$,
$(xf)(n)=f(nx)$,
$x\in \g_w$, $f\in U(\g_w)^*$,
$n\in U(\g_w)$.

If $V$ is a $\g$-module
and $w\in \eW$,
we obtain
a new $\g$-module
by twisting the action on $V$
as  $X\cdot v=w\inv(X)v$,
$X\in \g$.
The module obtained in this way
we shall  denote by
$
\phi_w(V)$.

Arkhipov \cite{Ark} defined a twisting functor
$T_w: \BGG_{\kappa}^{[\Lam]}\rightarrow
\BGG_{\kappa}^{[w\circ \Lam]}$,
$w\in W$,
%on the category of $\g$-modules
by
\begin{align*}
T_w(V)=S_w\*_{U(\g)}\phi_w( V).
\end{align*}
Let $w=s_{j_1}\dots s_{j_{\ell}}$
be a reduced expression of $w\in W$.
Then,
we have
\begin{align}\label{eq:Tw_dec}
T_w=T_{s_{j_1}}\circ  \dots \circ
T_{s_{j_{\ell}}}.
\end{align}

We extend the functor $T_w:\BGG_{\kappa}^{[\Lam]}\rightarrow
\BGG_{\kappa}^{[w\circ \Lam]}$
for $w\in \eW$
as follows:
For $x\in \eW_+$ and a $\g$-module $V$,
let $T_x(V)$
be the
$\g$-module obtained from $\phi_x(V)$
by twisting the action
as
$\Dg\cdot v=(\Dg +
\bra x(\rho)-\rho, \Dg\ket \id )v$ and
$X\cdot v=X v$
($X\in [\g,\g]$,
$v\in \phi_x(V)$).
Then,
$T_x(M(\lam))=M(x\circ \lam)$
($x\in \eW_+$).
Set for $\widetilde{w}=xw\in \eW$
($x\in \eW_+,w\in W$),
\begin{align}\label{eq:T_w_def_extend}
T_{\widetilde{w}}(V)=T_x\left(T_w(V)\right).
\end{align}
%Then,
%$V\rightsquigarrow T_w(V)$,
%$w\in \eW$
%defines a functor
%on the category of $\g$-modules.
Note that
\begin{align}\label{eq:T_w_as_g_w}
T_w(V)=U(\g_w)^* \*_{U(\g_w)}\phi_w(V)
\end{align}
as $U(\g_w)$-modules.
%We have:
%$T_w(V)\in \BGG_{\kappa}^{[w\circ \Lam]}$
%for $V\in \BGG_{\kappa}^{[\Lam]}$.

\subsection{Proof of
Theorem \ref{Th:cofreeness_of_Verma}
and Theorem \ref{Th:global_cofreeness}}
\begin{proof}[Proof of
Theorem \ref{Th:cofreeness_of_Verma}]
By
the proof of \cite[Proposition 6.3 (i)]{AL},
one sees that
\begin{align}\label{eq:verma-twist-verma}
M(\lam)=T_w\left(M(w\inv\circ \lam)\right)
~\text{ (if $\bra \lam+\rho,\alpha\che\ket
\not
\in \Z_{\geq 1}$ for all
$\alpha\in \prroots\cap w (\nrroots)$ )}
\end{align}
for $w\in \eW$.
Since $M(\lam)$, $\lam\in \dual{\h}$, is free
over
$w\inv (\g_w)\subset \g_-$,
Theorem
\ref{Th:cofreeness_of_Verma}
%the following Lemma
immediately follows from
\eqref{eq:T_w_as_g_w}
and \eqref{eq:verma-twist-verma}.
%:
%\begin{Lem}\label{Lem:cofreeness_of_Verma}
%Suppose that
%$\bra \lam+\rho,\alpha\che\ket\not\in \Z_{\geq 1}$
%for all
%$\alpha\in \prroots\cap w(\nrroots)$.
%Then,
%$M(\lam)$ is cofree as a
%$\g_w$-module.
%\end{Lem}
\end{proof}

\begin{proof}[Proof of Theorem
\ref{Th:global_cofreeness}]
It
is easy to see that
Proposition reduces the case when $w\in W$.
We shall proceed  by induction on $\ell(w)$
for $w\in W$.

The case when $\ell(w)=1$   follows from
\cite[Lemma 4.1]{KT2}.
Let
$w=ys_{\alpha}$,
$\alpha\in \Pi$,
$\ell(w)=\ell(y)+1$.
Set $\beta=y(\alpha)\in \prroots$.
Then,
\begin{align}
&\prroots\cap w(\nrroots)
=\prroots\cap y(\nrroots)
\sqcup \{ \beta\}, \label{eq:Lem;added1}\\
&\g_w^t=\g_y^t\+ \C x_{-\beta},\quad
[\g_w^t,\g_y^t]\subset \g_y^t.
\end{align}
Here,
$x_{-\beta}$ is a root vector of $\g$
of root $-\beta$.

Let
$\Lam$ be as in Theorem
and $V\in \BGG_{\kappa}^{[\Lam]}$.
By \eqref{eq:Lem;added1}
and the induction
hypothesis,
$V$ is free over $\g_y^t$.
We shall show that
$V/\g_y^t V$
is free over $\C x_{-\beta}$:

Let
$V'=\phi_y(T_{y\inv}(V))$.
Since
$T_{y\inv}(V)%=U(\g_{y\inv})^*\*_{U(\g_{y\inv})}\phi_{y\inv}(V)
=\phi_{y\inv}\left(U(\g_y^t)^*\*_{U(\g_y^t)}
V
\right)
$,
it follows that
$V'
=U(\g_y^t)^*\*_{U(\g_y^t)}
V$.
The freeness of
$V$ over $\g_y^t$
implies that
\begin{align}\label{eq:Lem-temp-iso-1}
(V')^{\g_y^t}
=\C 1^*\*_{U(\g_y^t)}
V
\cong V/
\g_y^tV
\end{align}
where
$(V')^{\g_y^t}
=H^0(\g_y^t,
V')\subset V'$
and
$1^*\in
U(\g_y^t)^*
$
is the dual element of $1\in U(\g_y^t)$.

We claim that
\eqref{eq:Lem-temp-iso-1}
is an isomorphism
of $\C x_{-\beta}$-modules.
Indeed,
one can show that
$x_{-\beta}$ acts on $V'
=U(\g_y^t)^*\*_{U(\g_y^t)}
V$ as
\begin{align}\label{eq:formula-action-inproof1}
x_{-\beta}(f\* v)=(\ad^*(x_{-\beta})f)\*v
+f\* x_{-\beta}v
\quad (f\in U(\g_y^t)^*,
~v\in V),
\end{align}
where $(\ad^*(x_{-\beta})f)(n)=-f([x_{-\beta},n])$.

Because
$T_{y\inv}(V)\in \BGG_{\kappa}^{[y\inv \circ \Lam]}$
and
$\bra y\inv \circ \Lam,\alpha\ket \not\in \Z$,
it follows that
$T_{y\inv}(V)$
is free over $\C x_{-\alpha}$ by \cite[Lemma 4.1]{KT2},
and thus $V'$ is free over $\C x_{-\beta}$.
Therefore,
$(V')^{\g_y^t}\subset
V'$ is also free over $\C x_{-\beta}$
since
$U
(\C x_{-\beta})=\C[x_{-\beta}]$
is a principal ideal domain.
Hence we conclude
that $V/
\g_y^tV$ is
free over $\C x_{-\beta}$,
proving
$H_i(\C x_{-\beta},V/\g_y^t V)=0$
for $i\ne 0$.
But then
the Hochschild-Serre spectral sequence
for
the ideal $\g_y^t\subset \g_w^t$
proves that
$H_i(\g_w^t,V)=0$
for $i\ne 0$.
This proves
Theorem
\ref{Th:cofreeness_of_Verma}
by
Proposition \ref{Pro:criterion-of-freeness}.
\end{proof}
\begin{Rem}\label{Rem:exact-twist}
Let $w\in \eW$
and
$\Lam\in \dual{\h}_{\kappa}$
as in Theorem \ref{Th:global_cofreeness}.
Then,
one can prove that
the functor
$T_w$ defines an equivalence
of
categories $\BGG_{\kappa}^{[w\inv \circ \Lam]}\isomap
\BGG_{\kappa}^{[\Lam]}$
such that
$T_w(M(w\inv\circ \lam))=M(\lam)$,
$T_w(L(w\inv \circ \lam))=L(\lam)$
for $\lam\in W^{\Lam}\circ \Lam$.
%see \cite[Theorem 2.1]{S}.
\end{Rem}
%By
%applying Proposition \ref{Pro:global_cofreeness}
%to $w=t_{\srho\che}$ and $w_0$,
%we obtain the following result
%which
% will be
%needed in section \ref{section:estimate}.
%\begin{Pro}\label{Co:freeness}
%Let $\Lam\in \dual{\h}_{\kappa}$.
%\begin{enumerate}
%\item
%Suppose that
%$\bra \Lam+\rho,\alpha\che\ket\not\in \Z$
%for all
%$\alpha\in \prroots\cap t_{\srho\che}(\nrroots)$.
%Then,
%any object in $\BGG^{[\Lam]}_{\kappa}$
%is free as a $t_{\srho\che}(\g_+)\cap \g_-$-module.
%\item
%Suppose that
%$\bra \Lam+\rho,\alpha\che\ket\not\in \Z$
%for all $\alpha\in \sproots$.
%Then,
%any object in $\BGG^{[\Lam]}_{\kappa}$
%is free as a $\snn$-module.
%\end{enumerate}
%\end{Pro}

\section{The BRST complex}
\label{section:complex}
In this section
we collect necessary information
from \cite{FB,FKW,FF_W,Feigin}
about the BRST
complex of the quantized Drinfeld-Sokolov reductions.
\subsection{Notations}
Let $\n$ denote $\Ln{+}
$ or $\Ln{-}$.
Here,
$\Ln{\pm}
=\sn_{\pm}\*\C[t,t\inv]\subset \g$
as in Introduction.

Let $\sn=\n\cap \sg$
and $\bar{\n}^t=\n^t\cap \sg$.
Then,
$\sg=\sn\+ \sh\+ \sn^t$.
Set
$\n_{\pm}=\n\cap \g_{\pm}$
and $\n_{\pm}^t=\n^t\cap \g_{\pm}$.
Then,
$\n=\n_-\+\n_+$
and $\n^t=\n^t_-\+ \n^t_+$.
%\begin{align*}
%&\n_-=\snp\*\C[t\inv]t\inv,
%\quad \n_+=\snp\* \C[t],\\
%&\n_-^t=\snn\* \C[t\inv],\quad
%\n_+^t=\snn\* \C[t]t
%\end{align*}
%%for $\n=\Ln{+}$,
%\begin{align*}
%&\n_-=\snn\*\C[t\inv]
%\quad \n_+=\snn\* \C[t]t,\\
%&\n_-^t=\snp\* \C[t\inv]t\inv,,\quad
%\n_+^t=\snp\* \C[t]
%\end{align*}
%for $\n=\Ln{-}$.

Let $\sroots(\n)=\begin{cases}
\sproots&\text{(for $\n=\Ln{+}$)},\\
\snroots&\text{(for $\n=\Ln{-}$).}
\end{cases}$
%for $\n=\Ln{+}$
%and $\sroots(\n)=\snroots$ for $\n=\Ln{-}$.
\ For an $\ad \h$-stable subspace $\m$ of
$\g$,
let
$\rroots(\m)=\{\alpha\in \rroots;
\g_{\alpha}\cap \m\ne \{0\}\}$,
where $\g_{\alpha}\subset \g$ is the root space of root
$\alpha$.
%Thus,
%\begin{align*}
%\rroots(\n_+)=\{\alpha+n\delta;
%\alpha\in \sproots,
%n\in \Z_{\geq 0}\},\quad
%\rroots(\n^t_+)=\{-\alpha+n\delta;
%\alpha\in \sproots,
%n\in \Z_{> 0}\}
%\end{align*}
%for $\n=\Ln{+}$,
%and
%\begin{align*}
%&
%\rroots(\n_+)=\{-\alpha+n\delta;
%\alpha\in \sproots,
%n\in \Z_{> 0}\},\quad
%\rroots(\n^t_+)=\{\alpha+n\delta;
%\alpha\in \sproots,
%n\in \Z_{\geq  0}\}
%\end{align*}
%for $\n=\Ln{-}$.
Then,
$\rroots=\rroots(\n)\sqcup \rroots(\n^t)$,
$\rroots_{\pm}=\rroots(\n_{\pm})\sqcup
\rroots(\n^t_{\pm})$.

Let $$\H=\sh\*\C[t,t\inv]\+\C K\+\C \Dg$$
be the {\em Heisenberg subalgebra
of $\g$}.
Let $\b=\n\+ \H\subset \g$.
Then,
$\g=\n^t\+ \b=\n\+ \b^t$.
Set
$\b_-=\b\cap\g_-$,
$\b_+=\b\cap (\h\+\g_+)$
so that $\b=\b_-\+\b_+$.
Similarly ,
let $\b^t_-=\b^t\cap\g_-$
and $\b_+^t=\b^t\cap (\h\+\g_+)$.
Then,
$\g_-=\n_-\+\b_-^t=\n^t_-\+ \b_-$
and $\h\+\g_+=\n_+\+\b_+^t=\n^t_+\+\b_+$
(see Table 1).

Let
\begin{align*}
B=U(\b)\*\Lam (\n),
\end{align*}
where
$\Lam (\n)$ is the Grassmann algebra
of $\n$.
We regard $B$ as a $\C$-algebra
containing
$U(\b)$ and $\Lam (\n)$
as its subalgebras
such that
$[X,\omega]=\ad(X)(\omega)\in \Lam (\n)$
for $X\in \b$
and $\omega\in \Lam(\n)$.
Then,
$B
=U(\b)\cdot
\Lam (\n)=\Lam (\n)\cdot U(\b)$.
Let
$ N=U(\n)\* \Lam(\n)
=U(\n)\cdot \Lam(\n)=\Lam(\n)
\cdot U(\n)\subset B$.
Similarly,
we define algebras
$B_{\pm}=U(\b_{\pm})\*U(\n_{\pm})
\subset B$,
$N_{\pm}=U(\n_{\pm})\*U(\n_{\pm})\subset N$,
$B^t=U(\b^t)\*\Lam(\n^t)$,
$N^t=U(\n^t)\* \Lam(\n^t)$,
$B^t_{\pm}=U(\b^t_{\pm})\*\Lam(\n^t_{\pm})$,
$N^t_{\pm}=U(\n_{\pm}^t)\* \Lam(\n_{\pm}^t)$,
$w(B)=U(w(\b))\*\Lam(w(\n))$
($w\in \eW$)
and so on.

Let
${}^t:B\rightarrow B^t
$
be the
algebra
anti-isomorphism
induced by the Chevalley anti-isomorphism of $\g$.
\begin{table}
\begin{center}
\begin{tabular}{c|ccccc}
$\n$   &$\sn$   &$\n_+$ &$\n_-$&
$\b$& $\n^t$
\\\hline
$L\sn_+ $ &$\sn_+$&$\sn_+\*\C [t]$  &
$\sn_+\*\C [t\inv]t\inv$&
$L\sn_+\+ \H$&$L\sn_-$\\
$L\sn_- $ &$\sn_-$&$\sn_-\*\C[t]t$  &
$\sn_-\*\C [t\inv]$&
$L\sn_-\+ \H$&$L\sn_+$
\end{tabular}
\end{center}
\caption[Notations]{Notations}
\label{table}
\end{table}
\subsection{The Clifford algebra}
Let
$\sI=\{1,2,\dots,\rank \sg\}$.
Choose a basis $\{J_a;a\in
\sI\sqcup
\sroots\}$  of $\sg$
such that
$J_{\alpha}\in \sg_{\alpha}$,
$(J_{\alpha},J_{-\alpha})=1$
and $(J_{\alpha})^t=J_{-\alpha}$
($\alpha\in \sroots$).
Let $c_{a,b}^c$ be the
structure constant with respect to
this basis.
Then,
$c_{\alpha,\beta}^{\gamma}
=-c_{-\alpha,-\beta}^{-\gamma}$
($\alpha,\beta,\gamma\in \sproots$).
In the sequel,
we identify
$\dual{\n}$ with $\n^t$ via $(~,~)$
(observe $(\Ln{\pm})^t=\Ln{\mp})$.

Let
$\Cl$ be the {\em Clifford algebra}
associated to
$\n\+\dual{\n}=\n\+ \n^t
$
and its natural symmetric bilinear form.
Denote by
$\psi_{\alpha}(n)$,
$\alpha\in \sroots$,
$n\in \Z$,
the
generators of $\Cl$
which correspond to
the elements
$J_{\alpha}(n)=(J_{-\alpha}(-n))^*$.
Then,
\begin{align*}
\{\psi_{\alpha}(m),{\psi}_{\beta}(n)\}=
\delta_{\alpha+\beta,0}\delta_{m+n,0}\quad
(\alpha,\beta\in \sroots, m,n\in \Z).
\end{align*}
Here,
$\{x,y\}=xy+yx$.
The algebra $\Cl$   contains
$\Lam(\n),\Lam(\n^t)
$ as its subalgebras and $\Cl=\Lam(\n^t)\*\Lam(\n)$
as  $\C$-vector spaces.

Let $\F(\n)$  be the irreducible
representation of $\Cl$ generated a vector $\1$
such that
$\psi_{\alpha}(n)\1=0$
($\alpha \in \sroots$,
$n\in \Z$,
$\alpha+n\delta\in \prroots$).
Then,
$\F(\n)=\Lam(\n^t_-)\*\Lam(\n_-)$
as  $\C$-vector spaces.
Let
\begin{align}
\F^p(\n)=
\sum\limits_{i-j=p}\Lam^i(\n^t_-)\*\Lam^j(\n_-)
\subset \F(\n)
\quad (p\in \Z).\label{eq:identification-of-F}
\end{align}
Then,
$\F(\n)=\sum_{p\in \Z}\F^p(\n)$.
By definition,
\begin{align*}
\text{$
\psi_{\alpha}(n)\F^i(\n)\subset
\F^{i- 1}(\n)$
and $
{\psi}_{-\alpha}(n)\F^i(\n)\subset\F^{i+ 1}(\n)
\quad$($\alpha\in \sroots(\n),n\in \Z$).}
\end{align*}
Obviously,
$\F(\n)\cong \F(\n^t)$
as $\Cl$-modules,
but  their gradings are opposite.

Let
${}^t:\Cl\rightarrow \Cl$
be the
algebra anti-isomorphism
defined by
$\psi_{\alpha}(n)\mapsto
\psi_{-\alpha}(-n)$
($\alpha\in \sroots$,
$n\in \Z$).
Then,
there is a unique none-degenerate
bilinear form
$\bra~,~\ket_{\F}:\F(\n)\times
\F(\n)\rightarrow \C$
such that $\bra\1,\1\ket_{\F}=1$,
$\bra\omega v,v'\ket_{\F}=\bra
v,\omega^t v'\ket_{\F}$,
$\omega\in\Cl,v\in \F(\n),
v'\in \F(\n)$.
It is none-degenerate on
$\F^i(\n)\times
\F^{i}(\n)$.

\subsection{The complex $C(\n,V)$}
For $V\in \BGG_{\kappa}$,
let
$$C(\n,V)=
V\* \F(\n)=\sum_{i\in \Z}
C^i(\n,V),
\quad \text{where $C^i(\n,V)=V\* \F^i(\n)$}.$$
Define
the operator
$d_{\n}^{\st}$
on $C(\n,V)$ by
\begin{align}
&d_{\n}^{\st}=\sum_{\alpha\in
\sroots(\n),n\in \Z}
J_{\alpha}(-n){\psi}_{-\alpha}(n)
-\frac{1}{2}\sum_{\ud{\alpha,\beta,\gamma\in
\sroots(\n)}{ k+l+m=0}}
c_{\alpha,\beta}^{\gamma}
{\psi}_{-\alpha}(k){\psi}_{-\beta}(l)
\psi_{\gamma}(m)
\end{align}
Here,
$J_{\alpha}(-n)$ acts on the first factor $V$
and $\psi_{\alpha}(n)$
acts on the second factor $\F(\n)$.
Then,
$(d_n^{\st})^2=0$,
$d_{\n}^{\st}C^i(\n,V)\subset
C^{i+1}(\n,V)$.
The cohomology
\begin{align*}
&\coho{\bullet}(\n,V)
=H^{\bullet}(C(\n,V),d_{\n}^{\st}).
\end{align*}
is called the %Feigin's
{\em semi-infinite cohomology  of
$\n$ with
coefficients in $V$ }(\cite{Feigin}).

Define $\chi_{\n}\in \dual{\n}\subset
\Cl$
by
\begin{align*}
&\chi_{\n}=\sum_{\alpha\in
\sPi}{\psi}_{-\alpha}(1)
\quad \text{(for $\n=\Ln{+}$)},\quad
\chi_{\n}=\sum_{\alpha\in
\bar{\Pi}}{\psi}_{\alpha}(0)
\quad \text{(for $\n=\Ln{-}$).}%\label{eq:chi_in_CL}.
\end{align*}
%Then,
%it defines a character
%$\chi_{\n}:U(\n)\rightarrow \C$.
%Let  $\C_{\chi_{\n}}$
%(resp.~$\C_{\chi_-}$)
%be the one-dimensional representation of $
%U(\n)$
%defined by $\chi_{\n}$.
($\chi_{\n}$ was denoted by $\chi_{\pm}$
in Introduction).
Let $d_{\n}=d_{\n}^{\st}+\chi_{\n}$.
Then,
$\chi_{\n}^2=0$,
$\{\chi_{\n},d_{\n}\}=0$
and
$\chi_{\n}C^i(\n,V)\subset
C^{i+1}(\n,V)$.
In particular,
$d_{\n}^2=0$.
Define
\begin{align}
\Wcoho^{\bullet}(\n,V)=H^{\bullet}
(C(\n,V),d_{\n}).
\end{align}
It is called the
{\em cohomology %$\Wcoho^{\bullet}(\n,V)$
of the BRST complex of the
quantized Drinfeld-Sokolov reduction
for $\n$
%(^^ ^^ $+$" and ^^ ^^ $-$" reduction)
associated to  $V$} (\cite{FF_W,FKW,FB}).

\subsection{The weight space decomposition}
The space
$C(\n,V)$ is naturally  a $\h$-module,
see
\eqref{eq:identification-of-F}.
Let $C(\n,V)=\bigoplus_{\lam
\in \dual{\h}}
C(\n,V)^{\lam}$
be its weight space decomposition.
Then,
$\dim C(\n,V)^{\lam}<\infty$
for all $\lam\in \dual{\h}$
and $V\in \BGG_{\kappa}$.
We have:
\begin{align}
&d_{\n}^{\st}C(\n,V)^{\lam}
\subset C(\n,V)^{\lam},
\label{eq:d-st-presevers-weight-spaces}\\
&
\chi_{\n}C(\n,V)^{\lam}
\subset \sum_{\alpha\in
\sPi}C(\n,V)^{\lam-\alpha+\delta}
\quad\text{(for $\n=\Ln{+}$)},
\label{eq:chi-not-preseve-weight-spaces-+}\\
&
\chi_{\n}C(\n,V)^{\lam}
\subset \sum_{\alpha\in \sPi}C
(\n,V)^{\lam+\alpha}\quad
\text{(for $\n=\Ln{-}$)}.
\label{eq:chi-not-preseve-weight-spaces--}
\end{align}
By \eqref{eq:d-st-presevers-weight-spaces},
the complex
$(C(\n,V),d_{\n}^{\st})$
is a direct sum of
finite-dimensional subcomplexes
$C(\n,V)^{\lam}$,
$\lam\in \dual{\h}$.
Therefore,
$\coho{\bullet}(\n,V)$,
$V\in \BGG_{\kappa}$, admits a weight space
decomposition:
$
\coho{\bullet}(\n,V)=\bigoplus_{\lam
\in \dual{\h}}\coho{\bullet}(\n,V)^{\lam}
$,
$\dim \coho{\bullet}(\n,V)^{\lam}
<\infty$
($\lam\in \dual{\h}$).

Note that $\bra~,~\ket_{\F}$ induces
a non-degenerate paring
$C(\n,V^*)\times C(\n,V)
\rightarrow \C$
which is non-degenerate
on
$C^{i}(\n,V^*)^{\lam}\times C^{i}(\n,V)^{\lam}
$,
$i\in \Z$,
$\lam\in \dual{\h}$.
Thus,
\begin{align}\label{eq:identification-dual}
C^{i}(\n,V^*)=
C^{i}(\n,V)^*
\end{align}
as  $\C$-vector spaces,
where ${}^*$ is defined in \eqref{eq:graded-dual}.
Let
\begin{align}
d_{\n}^t=d_{\n^t}^{\st}+\chi_{\n}^t
\in \End C(\n,V).
\label{eq:transpose-of-d}
\end{align}
Here,
$d_{\n^t}^{\st}$ acts on $C(\n,V)$
by the identification
$C(\n,V)=C(\n^t,V)$,
and
\begin{align*}
&\chi_{\n}^t=\sum_{\alpha\in \sPi}{\psi}_{\alpha}(-1)
\quad (\text{for $\n=\Ln{+}$)},
\quad \chi_{\n}^t=\sum_{\alpha\in
\sPi}{\psi}_{-\alpha}(0)
\quad(\text{for $\n=\Ln{-}$)}.
\end{align*}
Then,
$(d_{\n}^t)^2=0$,
$d_{\n}^t C^{i}(\n,V)
\subset C^{i-1}(\n,V)$,
and
\begin{align}\label{eq:diffrential-of-dual}
(d_{\n}f)(v)=f(d_{\n}^tv)
\quad
(f\in C(\n,V^*),
v\in C(\n,V)
)
\end{align}
under the identification
\eqref{eq:identification-dual}.
\subsection{The action of $B$ on
$C(\n,V)$}
In the sequel
we follow \cite{FKW}
for the definition
of the normal ordering $:~:$.
Thus,
$:\psi_{\alpha}(n)\psi_{\beta}(m):=
\begin{cases}
\psi_{\alpha}(n)\psi_{\beta}(m)&(\alpha+n\delta\in
\nrroots),\\
-\psi_{\beta}(m)\psi_{\alpha}(n)
&(\alpha+n\delta\in \prroots).
\end{cases}$
and so on.
We have:
\begin{align*}
%\label{eq:change-of-order-in-normal-ordering}
:\psi_{\alpha}(n)\psi_{\beta}(m):
=-:\psi_{\beta}(m)\psi_{\alpha}(n):
\quad (\alpha,\beta\in \sroots,n,m\in \Z)
\end{align*}
Let
\begin{align}
&\widehat{J}_a(n)=J_{a}(n)+\sum_\ud{\beta,\gamma\in
\sroots(\n)
}{ k\in \Z}c_{a,\beta}^{\gamma}
:\psi_{ \gamma}(n-k){\psi}_{-\beta}(k):~
\text{($a\in \sroots(\n)\sqcup \sI$
and $n\in \Z$)},
\\
&\Dtot=\Dg+\sum_{\alpha\in \sroots(\n),n\in
\Z}n:\psi_{\alpha}(n){\psi}_{-\alpha}(-n):.
\end{align}
Then,
for $V\in \BGG_{\kappa}$,
the correspondences
\begin{align}\label{eq:acion_of_Borel-+}
\begin{array}{ccccl}
\pi:& B=U(\b)\* \Lam (\n)
&\rightarrow
&
\End_{\C}C(\n,V)&\\
&J_{a}(n)&\mapsto &\widehat{J}_{a}(n)&
(a\in \sroots(\n)\sqcup \sI, n\in \Z)
\\
&K&\mapsto &\kappa \id&\\
&\Dg&\mapsto &\Dtot&\\
&\psi_{\alpha}(n)&\mapsto &\psi_{\alpha}(n)
&(\alpha\in \sroots(\n),
n\in \Z)
\end{array}
\end{align}
defines a representation of $
B$ on
$C(\n,V)$.
We have:
\begin{align}
&\widehat{J}_{\alpha}(n)=\{
d_{\n}^{\st},\psi_{\alpha}(n)\}
\quad \left(\text{$\alpha\in \sroots(\n)$,
$n\in
\Z$}\right)
\label{eq:def_of_J_hat},\\
&[d_{\n}^{\st},\widehat{J}_{a}(n)]=0
\quad
(a\in \sroots(\n)\sqcup \sI,n\in
\Z),\label{eq:commutativity-d-J}\\
&[\chi_{\n},\widehat{J}_{\alpha}(n)]=0
\quad (\alpha\in \sroots(\n),n\in \Z),
\label{eq:comu-rels-chi-J}\\
&C(\n,V)^{\lam}=
\{ v\in C(\n,V);
\pi(h)v=\bra \lam+h\che\Lam_0,h\ket v~
(h\in \h)\}.
\end{align}

Let $B^t$
act on $C(\n,V)$
via
the identification
$C(\n,V)=C(\n^t,V)$.
The representation of
$B^t$ obtained in this way we shall denote
by $\pi^t$.
Set
$\widehat{J}_{-\alpha}(n)=\pi^t(J_{-\alpha}(n))$
($\alpha\in \sroots(\n),
n\in \Z$).
Observe
$\pi_{|\H}=\pi^t_{|\H}$,
and under the identification
\eqref{eq:identification-dual},
\begin{align}\label{eq:transpose-hat-J}
(\pi(b)f)(v)=
f(\pi^t(b^t)v)\quad
(b\in B,
~f\in C(\n,V^*),
v\in C(\n,V)
).
\end{align}
We have:
\begin{align}
&[d_{\n}^{\st},\widehat{J}_{-\alpha}(n)]
=\sum_\ud{\beta\in \sroots(\n),
b\in \bar{I}\sqcup \sroots(\n^t)}{
k\in \Z}c_{\beta,
-\alpha}^b:{\psi}_{-\beta}(k)\widehat{J}_b(n-k):
-nk_{\alpha}{\psi}_{-\alpha}(n)\label{eq:com_re_d_J}
\end{align}
for $\alpha\in \sroots
(\n),n\in \Z$,
where $k_{\alpha}=\kappa-
h\che-\sum\limits_{\beta,\gamma
\in \sroots(\n)}c_{\alpha,\beta}^{\gamma}
c_{-\alpha,-\beta}^{-\gamma}\in \C
$,
and
\begin{align}
&[\chi_{\n},\widehat{J}_a(n)]=
\sum_{\beta,\gamma\in \sroots(\n)}c_{a,\beta}^{\gamma}
\chi_{\n}(J_{\gamma}(-1)){\psi}_{-\beta}
(n+1)\quad \text{(for $\n=\Ln{+}$)}
\label{eq:com-rel-chi-J-+}\\
&[\chi_{\n},\widehat{J}_a(n)]=
\sum_{\beta,\gamma\in \sroots(\n)}c_{a,\beta}^{\gamma}
\chi_{\n}(J_{\gamma}(0)){\psi}_{-\beta}
(n)\quad \text{(for $\n=\Ln{-}$)}
\label{eq:com-rel-chi-J--}
\end{align}
for $a\in \sroots(\n^t)\sqcup \sI$
and $n\in \Z$.

\subsection{The degree operator
$\DW_{\n}$}
Define
\begin{align*}
&\DW_{\n}=\Dtot +\pi({\srho\che})\quad
\quad\text{(for $\n=\Ln{+}$),}\\
& \DW_{\n}=\Dtot
+\left(\frac{1}{2}|\srho\che|^2\kappa
-\bra \srho,\srho\che\ket\right)\id
\quad\text{(for $\n=\Ln{-}$)}.
\end{align*}
Set
\begin{align*}
&\Whw{\n}{\lam}=\bra \lam, \bar{\rho}\che +\Dg\ket
~\text{(for $\n=\Ln{+})$, }\quad
\Whw{\n}{\lam}=\bra
t_{\srho\che}\circ \lam,
\bar{\rho}\che +\Dg
\ket
~ \text{(for $\n=\Ln{-})$.}
\end{align*}
Then,
${\DW_{\n}}$
acts as the multiplication
by $\Whw{\n}{\lam}\in \C$
on the weight space
$C(\n,V)^{\lam}$,
$\lam\in \dual{\h}_{\kappa}$.
Let
$C(\n,V)_a=
\{v\in C(\n,V);
\DW_{\n}v=av\}$
for $a\in \C$.
Clearly,
$
C(\n,V)=
\bigoplus_{a\in\C}C
(\n,V)_{a}
$
and \begin{align}
&C(\n,V)_a=\bigoplus_\ud{\lam\in
\dual{\h}}{
\Whw{\n}{\lam}=a}C(\n,V)^{\lam}.
\label{eq:description_of-eigen-values}
\end{align}
By \eqref{eq:d-st-presevers-weight-spaces},
\eqref{eq:chi-not-preseve-weight-spaces-+}
and \eqref{eq:chi-not-preseve-weight-spaces--},
it follows that
$d_{\n}^{\st}C(\n,V)_a
\subset
C(\n,V)_a$
and $\chi_{\n}C(\n,V)_a\subset
C(\n,V)_a$.
Therefore,
\begin{align}
\text{$\Wcoho^{\bullet}(\n,V)=
\bigoplus_{a\in \C}
\Wcoho^{\bullet}(\n,V)_a,$
}
\end{align}
where
$\Wcoho^{\bullet}(\n,V)_a
=H^{\bullet}(C(\n,V)_a,d_{\n})$.
\begin{Rem}\label{Rem:Weyl-grop-action}
\begin{enumerate}
\item
The operator  $\DW_{\n}$
is the semisimplification
of
$-L_0^{\pm}$ defined in
\cite[3.1]{FKW} up to constant shift.
\item
The eigenspace
$C(\n,V)_a$
of $\DW_{\n}$ is not
necessarily
finite-dimensional.
\end{enumerate}
\end{Rem}

\subsection{The Weyl group
action}\label{subsection:Wey-group-action}
The group $\eW$ acts naturally
on $\Cl$.
Let $\eW$
act on $U(\g)\*\Cl$
by $w(u\* \omega)=
w(u)\* w(\omega)$.

Let
$w\in \eW$.
If  $V$ is a $U(\g)\*\Cl$-module,
we obtain
a new $U(\g)\*\Cl$-module
by twisting its action on $V$
as  $A\cdot v=w\inv(A)v$.
The module obtained in this way
we shall also denote by
$
\phi_w(V)$.
Then,
the action
$\pi_{|\H}$ of $\H$ is well-defined on
$\phi_w(C(\n,V))$.
By direct calculation,
one gets the following proposition.
\begin{Pro}\label{Pro:tristing-heisenbeg}
For $w\in \eW$
and $v\in C(\n,V)$,
\begin{align*}
&\pi(h)\cdot\phi_w( v)=\phi_w\left(\pi(w\inv(h))v
+\bra \rho,w\inv(h)-h\ket v\right)
\quad (h\in \h),\\
&\widehat{h}(n)\cdot \phi_w(v)=
\phi_w(\widehat{h}(n)v)\quad
(h\in \sh,n\in \Z\backslash \{0\})
\end{align*}
where
$\phi_w(v)$ denotes
the image of $v$
in $\phi_w(C(\n,V))$.
\end{Pro}
\section{Cohomology
associated to Verma modules
}
\label{section:vanishing_of_projectives}
In this section we review
the results obtained in \cite[14.2]{FB}
for our case.

\subsection{The decomposition
of $C(\n,M(\lam))$}
Fix  $\lam\in \dual{\h}$.
Let %$\lam\in \dual{\h}$
%and
$v_{\lam}$  be the highest weight vector of
$M(\lam)$.
Let $
|\lam \ket=v_{\lam}\*\1\in C
(\n,M(\lam))
$.
Then,
\begin{align}\label{eq:hightest}
d_{\n}^{\st}|\lam\ket=0,\quad
\chi_{\n}|\lam\ket=0.
\end{align}
\begin{Pro}\label{Pro:dec-of-M(lam)}
The map defines by
\begin{align}\label{eq:dec-of-M(lam)}
\begin{array}
{ccc}
N_-
\* B^t_-
&\longrightarrow &C(\n,M(\lam))\\
n\* b&\longmapsto &
\pi(n)\cdot \pi^t(b)|\lam\ket
\end{array}
\end{align}
gives an isomorphism
of $\C$-vector spaces.
\end{Pro}
\begin{proof}
Observe
that $M(\lam)=U(\n_-)U(\b^t_-)v_{\lam}$.
Thus,
by comparing the dimension
of weight spaces of both sides,
it follows that
it is sufficient to show that
\eqref{eq:dec-of-M(lam)} is surjection.
By definition,
\begin{align*}
\pi(N_-)\pi(B_-)|\lam\ket=U(\n_-)\cdot
\Lam(\n_-)\cdot U(\b^t_-)\cdot
\Lam(\n_-^t) |\lam\ket.
\end{align*}
But we have
$U(\b^t_-)\cdot
\Lam(\n_-)=\Lam(\n_-)\cdot U(\b^t_-)$.
This can be seen
by the commutation relations
\begin{align*}
[\widehat{J}_a(z),\psi_{\beta}(z)]=
-\sum_{\beta,\gamma\in
\sroots(\n)}
c_{a,-\gamma}^{-\beta}\psi_{\gamma}(n+m)
\quad (a\in \sroots(\n^t)\sqcup \sI,
\beta\in \sroots(\n),
n,m\in \Z).
\end{align*}
On the other hand,
$\Lam(\n_-)\cdot \Lam(\n_-^t)|\lam\ket
=\left(\C v_{\lam}\right)\* \F(\n)$.
Thus,
it is enough to show that
\begin{align*}
U(\n_-)\cdot
U(\b^t_-)\cdot
\left(\left(\C v_{\lam}\right)\* \F(\n)\right)
=M(\lam)\* \F(\n).
\end{align*}
But
$U(\b^t_-) v_{\lam}$ is a free
$\b^t_-$-submodule of $M(\lam)$.
Thus,
$
U(\b^t_-)\cdot
\left(\left(\C v_{\lam}\right)\* \F(\n)\right)=
\left(U(\b^t_-) v_{\lam}\right)\* \F(\n)
$.
Similarly,
we get
$
U(\n_-)\cdot \left(U(\b^t_-) v_{\lam}\right)\* \F(\n)
=M(\lam)\* \F(\n)
$.
\end{proof}
\subsection{
The subcomplex ${C(\n,\lam)_0}$}
\label{subsection:C_0}
We define
the subspace
$C(\lam)_0={C(\n,\lam)_0}$
of $C(\n,M(\lam))$,
$\lam\in \dual{\h}_{\kappa}$,
by
\begin{align}
{C(\lam)_0}
&=\pi^t(B^t) |\lam\ket
\end{align}
Since
$\pi^t$ defines one-dimensional representation
of $B^t_+$
on $\C |\lam\ket$,
it follows that
\begin{align*}
C(\lam)_0=\pi^t(B^t_-)
|\lam\ket
= U(\b^t_-)\cdot \Lam (\n^t_-)\cdot|\lam\ket=\Lam
(\n^t_-)\cdot U(\b^t_-)\cdot  |\lam\ket.
\end{align*}
\begin{Pro}\label{Pro:iso-map-decomposition}
$
{C(\lam)_0}=
B^t\*_{
B^t_+}\C |\lam\ket
$.
\end{Pro}
\begin{proof}
It is sufficient to
show
that
the multiplication
map
$B^t_-
\* \C|\lam\ket
\rightarrow
\pi^t(B^t_-)
|\lam\ket={C(\lam)_0}$
is an injection.
But this easily follows form
Proposition \ref{Pro:dec-of-M(lam)}.
\end{proof}
By \cite[14.2]{FB},
${C(\lam)_0}$
is a subcomplex of $C(\n,M(\lam))$.
This can be seen
from \eqref{eq:hightest}
and the commutation relations
\eqref{eq:com_re_d_J},
\eqref{eq:com-rel-chi-J-+},
\eqref{eq:com-rel-chi-J--}.

Let
$\left({C(\lam)_0}\right)_a
={C(\lam)_0}\cap
C\left(\n,M(\lam)\right)_a$,
$a\in \C$.
The following Proposition is easy to
see.
\begin{Pro}\label{Pro:finite-dimensioanl}
For $\lam\in \dual{\h}$,
$C(\lam)_0=\bigoplus\limits_{a\in\Whw{\n}{\lam}
+\Z_{\leq 0}}
\left({C(\lam)_0}\right)_a$
and
$\left({C(\lam)_0}\right)_a$
is finite-dimensional for all $a$.
\end{Pro}

Define
the subspace
$C(\n_-)'$
of $C(\n,M(\lam))$
by
\begin{align*}
C(\n_-)'=
\pi(N_-) |\lam\ket
=U(\n_{-})\cdot  \Lam (\n_{-})\cdot |\lam\ket
=\Lam (\n_{-})\cdot U(\n_{-})\cdot |\lam\ket.
\end{align*}
Then,
$C(\n_-)'\cong N_-
$ as $\C$-vector spaces.
This can be seen in the same way as
Proposition \ref{Pro:iso-map-decomposition}.
By \eqref{eq:hightest},
\eqref{eq:def_of_J_hat},
\eqref{eq:commutativity-d-J},
\eqref{eq:comu-rels-chi-J}
and the fact
\begin{align}\label{eq:com-re;chi-J-+}
\{\chi_{\n},\psi_{\alpha}(n)\}=\chi_{\n}(J_{\alpha}(n))\quad
(\alpha\in \sroots(\n),n\in \Z),
\end{align}
it follows that
$C(\n_-)'
$
is a subcomplex of  $C(\n,M(\lam))$.
It is easy to see that
this complex
does not depend on $\lam\in \dual{\h}_{\kappa}$.
\begin{Pro}\label{Pro:dec-of-complex-Verma}
\cite[14.2]{FB}
The map in Proposition \ref{Pro:dec-of-M(lam)}
defines an isomorphism
\begin{align*}
C(\n,M(\lam))\isomap
C(\n_-)'\* {C(\lam)_0}
\end{align*}
of complexes.
\end{Pro}

Though the following Proposition
is proved
in \cite[14.2]{FB}
in slightly different setting,
the same proof applies.
\begin{Pro}[{\cite[14.2]{FB}}]
\label{Pro:decompsotion-FB}
\begin{enumerate}
\item
$H^{i}\left({C(\lam)_0}
\right)=\{0\}$
for $i\ne 0$.
\item $H^i(C(\n_-)')=
H^i(C(\n_-)', d_{\n}^{\st})
=\begin{cases}
\C &(i=0)\\
\{0\}&(i\ne 0).
\end{cases}
$
\end{enumerate}\end{Pro}
\begin{Rem}
In the proof of Proposition \ref{Pro:decompsotion-FB}
(1),
one uses the fact that
\begin{align*}
&P(C(\lam)_0)\subset
\{ \mu\in \dual{\h}_{\kappa};
\bra \lam-\mu,\srho\che\ket
\in\Z_{\geq 0}\}
\quad \text{(for $\n=\Ln{+}$)},\\
&P(C(\lam)_0)\subset
\{ \mu\in \dual{\h}_{\kappa};
\bra \mu-\lam,\srho\che\ket\in\Z_{\geq 0}\}
\quad \text{(for $\n=\Ln{-}$)}
\end{align*}
to assure the convergency
of the spectral sequence described in \cite[14.2.8]{FB}.
\end{Rem}
Proposition \ref{Pro:dec-of-complex-Verma}
and
Proposition \ref{Pro:decompsotion-FB}
imply:
\begin{Th}\label{Th:vanishing.for.Verma}
For $\lam\in \dual{\h}$,
$\Wcoho^{i}\left(\Ln{\pm},M(\lam)\right)=\{0\}$
$(i\ne 0)$.
\end{Th}
\begin{Rem}\label{Rem:character-of-VErma}
Set
$\ch
\Wcoho^{i}\left(\n,V\right)=
\tr_{\Wcoho^{i}\left(\n,V\right)}
q^{\DW_{\n}}
$
when
$\dim \Wcoho^{i}\left(\n,V\right)_a<\infty$
for all $a$.
It is easy to see that
$
\ch \Wcoho^{0}\left(\n,M(\lam)\right)
=\frac{q^{\Whw{\n}{\lam}}}{
\prod_{i\geq 1}(1-q^i)^{\rank \sg}}
$
for $\lam\in \dual{\h}$,
see {\cite[14.2]{FB}}.
\end{Rem}
\section{Cohomology associated to
duals of Verma modules}
\label{section:vanishing_of_dual_of_Verma}
In this section
we prove the
vanishing of
$\Wcoho^{i}(\n,M(\lam)^*)$
for $i\ne 0$ under the certain
restriction
of $\lam$.
\subsection{Relative complex}
For $\lam\in \dual{\h}$,
define the subspace
$C(\lam)^0=C(\n,\lam)^0$ of
$C(\n,M(\lam)^*)$ by
\begin{align*}
&C(\lam)^0
= \left\{v\in C\left(\n,M(\lam)^*\right);
\begin{array}{l}
\widehat{J}_{\alpha}(n)v=\psi_{\alpha}(n)v=0\\
(\alpha\in \sroots(\n),n\in \Z,
\alpha+n\delta \in \rroots(\n_+))
\end{array}\right\}.
\end{align*}
It is the relative complex
(with respect to the differential $d_{\n}^{\st}$)
considered in \cite{FGZ}.
By \eqref{eq:def_of_J_hat},
\eqref{eq:commutativity-d-J}
and \eqref{eq:com-re;chi-J-+},
it follows that
$d_{\n}C(\lam)^0
\subset C(\lam)^0$.
Thus, it
is a subcomplex of $C(\n,M(\lam)^*)$.
\begin{Pro}\label{Pro:identification-of-C-top-0}
For $\lam\in \dual{\h}$,
$
C(\lam)^0=
\Big(M(\lam)^*\*
\Lam (\n_-)\Big)^{ \n_+}
\subset C(\n,M(\lam)^*)
$
and the restriction of
$d_{\n}$ on $C(\lam)^0$
is $
d_{\n_-}=d_{\n_-}^{\st}+\chi_{\n}$.
Here,
$\n_+$ acts on $\Lam (\n_-)$
via the identification
$\n_-=\n/\n_+$,
and $d_{\n_-}^{\st}$
is the differential
of $\n_-$-homology,
that is,
\begin{align*}
d_{\n_-}^{\st}&=\sum_{\alpha+n\delta
\in \rroots(\n_-)}
J_{\alpha}(n)\psi_{-\alpha}(-n)\\ &\quad\quad\quad-
\frac{1}{2}\sum_{\alpha+ k\delta,\beta+l\delta
\in \rroots(\n_-)}c_{\alpha,\beta}^{\gamma}
\psi_{-\alpha}(-k)\psi_{-\beta}(-l)\psi_{\gamma}(k+l).
\end{align*}
\end{Pro}
\begin{proof}
Clearly,
$C(\lam)^0$
is contained
in the subspace
\begin{align*}
M(\lam)^*\*
\Lam(\n_-)
=\left\{v\in
C(\n,M(\lam)^*);
\psi_{\alpha}(n)v=0
\quad (\alpha\in \sroots(\n),n\in \Z,
\alpha+n\delta\in \prroots)
\right\}.
\end{align*}
It is easy to see that
the operators
$\widehat{J}_{\alpha}(n)$
($ \alpha\in \sroots(\n),n\in \Z,
\alpha+n\delta\in \prroots$)
preserve this subspace
and their action
coincide with the one
via the identification
$\n_-=\n/\n_+$.
Hence,
it follows that
$C(\lam)^0=
\Big(M(\lam)^*\*
\Lam (\n_-)\Big)^{ \n_+}$
as  $\C$-vector spaces.
But then,
by the proof of
\cite[Theorem 2.2]{V1},
it follows that
the restriction
of $d_{\n}^{\st}$
to this subspace is $d_{\n_-}^{\st}$.
\end{proof}

Define a subspace
$C(\lam)^t_0=C
(\n,\lam)^t$
of
$C(\n,M(\lam))$
by \begin{align*}
{{C(\lam)^t_0}}
=B\*_{
B_+} \C |\lam\ket,
\end{align*}
see Proposition
\ref{Pro:iso-map-decomposition}.
Then,
$\chi_{\n}^tC(\lam)^t_0
\subset C(\lam)^t_0$.
We view
${{C(\lam)^t_0}}$
as a complex with differential $d_{\n}^t$,
where
$d_{\n}^t$ is defined in
\eqref{eq:transpose-of-d}.
%Here,
%the fact that
%$d_{\n}^tC(\lam)^t_0
%\subset C(\lam)^t_0$,
%is easily seen.
\begin{Pro}\label{Pro:identification-dual-0}
For $\lam\in \dual{\h}$,
$
C(\lam)^0={\left( C(\lam)^t_0\right)}^*
$ as a complex.
\end{Pro}
\begin{proof}
follows from
\eqref{eq:identification-dual},
\eqref{eq:transpose-hat-J} and
Proposition
\ref{Pro:dec-of-M(lam)}.
\end{proof}
\begin{Pro}\label{Pro:reduction-of-dual}
For $\lam\in \dual{\h}$,
\begin{align*}
\Wcoho^{\bullet}\left(\n,M(\lam)^*\right)
=H^{\bullet}(C(\lam)^0)
\left(=H^{\bullet}\left(
{\left( C(\lam)^t_0\right)}^*\right)\right).
\end{align*}
\end{Pro}
\begin{proof}
The proof can be done using
the corresponding statement to
Proposition \ref{Pro:dec-of-complex-Verma}.
Or one can apply \cite[Theorem 2.2]{V1}.
Indeed,
by Proposition \ref{Pro:identification-of-C-top-0},
the complex
$C(\lam)^0$ is nothing but the $E_1^{\bullet,0}$-row
of the Hochschild-Serre spectral sequence
for $\n_+\subset \n$
in
\cite[Theorem 2.2]{V1}.
But since $M(\lam)^*$ is a cofree $\n_+$-module,
it follows that
$E_1^{\bullet,q}=0$ for $q\ne 0$.
Thus this spectral sequence collapses at $E_2=
H^{\bullet}(C(\lam)^0)
=E_{\infty}$.
\end{proof}

We have:
\begin{align*}
(\DW_{\n}f)(v)=f(\DW_{\n}v)
\quad (f\in C(\lam)^0
={\left( C(\lam)^t_0\right)}^*,
v\in {C(\lam)^t_0}).
\end{align*}
Let
$
({{C(\lam)^t_0}})_{a}=
\{v\in {{C(\lam)^t_0}};
\DW_{\n}v=av\}
$.
Then,
$d_{\n}^t({{C(\lam)^t_0}})_{a}
\subset ({{C(\lam)^t_0}})_{a}$
and
$
{{C(\lam)^t_0}}=\bigoplus\limits_{a\in
\Whw{\n}{\lam}+\Z}({{C(\lam)^t_0}})_{a}
$.
Observe that
%on the contrary to
%Proposition \ref{Pro:finite-dimensioanl},
the eigenspace
$({{C(\lam)^t_0}})_{a}$ is not
finite-dimensional in general
(compare Proposition \ref{Pro:finite-dimensioanl}).
Below we shall define
a subspace
${\N\left(\lam\right)_0^t}
={\N\left(\n,\lam\right)_0^t}
\subset {{C(\lam)^t_0}}$
so that
the quotient ${C(\lam)^t_0}/
{\N\left(\lam\right)_0^t}
$ is a direct sum of finite-dimensional
eigenspaces of
$\DW_{\n}$.
The definition is different
for $\n=\Ln{+}$ and $\n=\Ln{-}$.

\subsection{The subspace
${\N\left(\lam\right)_0^t}$
for $\n=\Ln{+}$}
Let $\n=\Ln{+}$.
Observe that
\begin{align*}
&\n\supset
t_{\srho\che}(\n_+)\supset\n_+,
\quad\b\supset
t_{\srho\che}(\b_+)=\H_+\+t_{\srho\che}(\n_+)\supset
\b_+,
\\
&t_{\srho\che}(\b_+)\cap
\b_-=t_{\srho\che}(\n_+)\cap \n_-
=\haru \{
J_{\alpha}(n);\alpha
\in \sproots,-\height  \alpha\leq n\leq
-1\}.\end{align*}
Therefore,
we have inclusions of  algebras
$
B\supset
t_{\srho\che}(B_+)
\supset B_+
$.
Notice that
$\chi_n^t\in t_{\srho\che}(B_+)$.
Define
\begin{align*}
{\bar{C}(\lam)^t}
=t_{\srho\che}(B_+)
\*_{B_+}\C |\lam\ket.
\end{align*}
It is
a subspace of
${{C(\lam)^t_0}}$
spanned by  the vectors of the form
\begin{align*}
{\psi}_{\alpha_{r_1}}(-m_{1})
\dots {\psi}_{\alpha_{r_p}}(-m_{p})
\widehat{J}_{\alpha_{s_{1}}}(-n_{1})
\dots
\widehat{J}_{\alpha_{s_p}}(-n_{s})|\lam\ket
\end{align*}
with
$\alpha_{r_i}, \alpha_{s_i}\in \sproots$,
$1\leq m_{i}\leq \height \alpha_{r_i}$,
$1\leq n_{i}\leq \height \alpha_{s_i}$.
By definition,
$\chi_{\n}^t{\bar{C}(\lam)^t}
\subset {\bar{C}(\lam)^t}$
and
\begin{align}
&{\bar{C}(\lam)^t}
=\bigoplus\limits_{a\in \Whw{\n}{\lam}+\Z_{\geq 0}}
{\bar{C}(\lam)^t_{a}},
\quad \text{where ${\bar{C}(\lam)^t_{a}}
={\bar{C}(\lam)^t}\cap ({{C(\lam)^t_0}})_{a}$.}
\label{eq:d-weights-of-bar-C-+}
\end{align}
Define
the subspace
${\bar{\N}(\lam)^t}$
of ${\bar{C}(\lam)^t}$
by
\begin{align*}
&{\bar{\N}(\lam)^t}=
\sum_{
\bra \mu-\lam,\srho\che\ket>
0}\left(\bar{C}(\lam)^t\right)^{\mu}.
\end{align*}
Then,
$t_{\srho\che}(B_+)
\cdot
{\bar{\N}(\lam)^t}
\subset
{\bar{\N}(\lam)^t}$.
In particular,
$\chi_{\n}^t{\bar{\N}(\lam)^t}\subset {\bar{\N}(\lam)^t}$.
Define
\begin{align}
{{\N\left(\lam\right)_0^t}}=
B\*_{t_{\srho\che}(B_+)
}{\bar{\N}(\lam)^t}
\subset
{{C(\lam)^t_0}}.
\label{eq:def-of-N-+}
\end{align}
Then,
\begin{align}
\chi_{\n}^t{{\N\left(\lam\right)_0^t}}
\subset {{\N\left(\lam\right)_0^t}}.
\label{eq:chi-n-sub-N-+}
\end{align}
Observe that
${\bar{C}(\lam)^t}/{\bar{\N}(\lam)^t}$
is spanned by the image $\overline{|\lam\ket}$
of $|\lam\ket$.
We have
\begin{align}\label{eq:one-dim-+}
&d_{\n}^t\overline{|\lam\ket}=0,
\\& \widehat{J}_i(n)v=0~~(i\in \sI,n>0),\quad
\widehat{J}_{\alpha}(n)\overline{|\lam\ket}
={\psi}_{\alpha}(n)\overline{|\lam\ket}
=0~~(\alpha\in \sproots,n\geq -\height \alpha),
\\
&{{C(\lam)^t_0}}/{{\N\left(\lam\right)_0^t}}=
B\*_{
t_{\srho\che}(B_+)
}\C\overline{|\lam\ket}.
\label{eq:C/N-+}
\end{align}

\subsection{The subspace
${\N\left(\lam\right)_0^t}$
for $\n=\Ln{-}$}
Let $\n=\Ln{-}$
and $w_0$ be the longest element of $\sW$.
Then,
\begin{align*}
\n\supset w_0(\n^t_+)\supset \n_+
,
\quad \b\supset w_0(\b^t_+)
\supset \b_+.
\end{align*}
Thus,
$B\supset w_0(B_+^t)\supset B_+$.
Notice $\chi_{\n}^t\in w_0(B_+)$.
Define
\begin{align*}
{\bar{C}(\lam)^t}
=w_0(B^t_+)\*_{B_+}\C |\lam\ket
\subset {{C(\lam)^t_0}}.
\end{align*}
It is the span of the vectors of the form
\begin{align*}
\psi_{-\alpha_{r_1}}(0)
\dots \psi_{-\alpha_{r_p}}(0)
\widehat{J}_{-\alpha_{s_{1}}}(0)
\dots
\widehat{J}_{-\alpha_{s_p}}(0)|\lam\ket
\end{align*}
with
$\alpha_{r_i},\alpha_{s_i}\in \sproots$.
We have:
\begin{align}
&{{C(\lam)^t_0}}=B
\*_{w_0(B^t_+))
}{\bar{C}(\lam)^t},\label{eq:C-is-induce-from-bar--}
\end{align}
$\chi_{\n}^t{\bar{C}(\lam)^t}
\subset {\bar{C}(\lam)^t}$
and
\begin{align}
\label{eq:d-weights-of-bar-C--}
{\bar{C}(\lam)^t}
=\bar{C}(\lam)^t_{\Whw{\n}{\lam}}.
\end{align}

Define
the subspace
${\bar{\N}(\lam)^t}$
of ${\bar{C}(\lam)^t}$
by
\begin{align*}
{\bar{\N}(\lam)^t}=
\sum_{
\bra \lam-\mu,\srho\che\ket>
0}\left(\bar{C}(\lam)^t\right)^{\mu}
\end{align*}
Then,
$w_0(B^t_+)\cdot
{\bar{\N}(\lam)^t}\subset {\bar{\N}(\lam)^t}$,
in particular,
$\chi_{\n}^t{\bar{\N}(\lam)^t}
\subset {\bar{\N}(\lam)^t}$.
Define
\begin{align}
&{{\N\left(\lam\right)_0^t}}=
{B}\*_{w_0(B^t_+))
}{\bar{\N}(\lam)^t}\subset
{{C(\lam)^t_0}}.\label{eq:def-of-N--}
\end{align}
We have
\begin{align}
&\chi_{\n}^t
{{\N\left(\lam\right)_0^t}}
\subset {{\N\left(\lam\right)_0^t}},
\label{eq:chiN-N}\\
&{{C(\lam)^t_0}}/
{{\N\left(\lam\right)_0^t}}=
{B}\*_{w_0(B^t_+))
}\C\overline{|\lam\ket}.\label{eq:C/N--}
\end{align}
Here,
$\overline{|\lam\ket}$
is the image of $|\lam\ket$ in
${\bar{C}(\lam)^t}/{\bar{\N}(\lam)^t}$:
\begin{align*}
&d_{\n}^t\overline{|\lam\ket}=0,
\quad
\widehat{J}_i(n)v=0\quad(i\in \sI,n>0),
\\
&
\widehat{J}_{-\alpha}(0)\overline{|\lam\ket}
={\psi}_{-\alpha}(0)\overline{|\lam\ket}
=0\quad(\alpha\in \snroots).
\end{align*}

\subsection{${\N\left(\lam\right)_0^t}$
is a null subcomplex}
Let
$\left({{C(\lam)^t_0}}/
{{\N\left(\lam\right)_0^t}}\right)_{a}$ be the
$\DW_{\n}$-eigenspace of
${{C(\lam)^t_0}}/
{{\N\left(\lam\right)_0^t}}$
of eigenvalue $a\in
\C$.
The following proposition
is easy to see
by \eqref{eq:C/N-+}
and \eqref{eq:C/N--}.
\begin{Pro}\label{Pro:finite-dim-of-C/N}
For $\lam\in \dual{\h}$,
$C(\lam)^t_0/
\N\left(\lam\right)_0^t
=\bigoplus\limits_{a\in \Whw{\n}{\lam}+\Z_{\leq 0}}
\left({{C(\lam)^t_0}}/
{{\N\left(\lam\right)_0^t}}\right)_a$
and
$\left({{C(\lam)^t_0}}/
{{\N\left(\lam\right)_0^t}}\right)_{a}$
is finite-dimensional for all $a\in \C$.
\end{Pro}

The proof of the following proposition
will be given
in
\ref{subsection:The
spectral seqeunce I}
and \ref{subsection:spectral-I}.
\begin{Pro}\label{Pro:later}
{\rm (1)}
Let $\n=\Ln{+}$
For $\lam\in \dual{\h}$,
$d_{\n^t}^{\st}{{\N\left(\lam\right)_0^t}}
\subset {{\N\left(\lam\right)_0^t}}$.
Moreover,
if
$\bra \lam+\rho,\alpha\che\ket\not\in \Z_{\geq 1}$
for
all $\alpha\in \prroots\cap t_{\srho\che}(\nrroots)$,
then
$H_{\bullet}
\left({{\N\left(\lam\right)_0^t}},d_{\n^t}^{\st}
\right)\equiv 0$.
{\rm (2)}
Let $\n=\Ln{-}$.
For $\lam\in \dual{\h}$,
$d_{\n^t}^{\st}{{\N\left(\lam\right)_0^t}}
\subset {{\N\left(\lam\right)_0^t}}$.
Moreover,
if
$\bra \lam+\rho,\alpha\che\ket\not\in
\Z_{\geq 1}$
for all  $\alpha\in \sproots$,
then
$H_{\bullet}
\left({{\N\left(\lam\right)_0^t}},d_{\n^t}^{\st}
\right)\equiv 0$.
\end{Pro}
By
\eqref{eq:chi-n-sub-N-+},
\eqref{eq:chiN-N}
as Proposition  \ref{Pro:later},
we have an exact sequence
$
0\rightarrow
{{\N\left(\lam\right)_0^t}}
\rightarrow
{{C(\lam)^t_0}}
\rightarrow
{{C(\lam)^t_0}}/{{\N\left(\lam\right)_0^t}}
\rightarrow 0
$
of complexes.
Therefore,
we get
the following  exact sequence of complexes:
\begin{align}\label{eq:exact}
0\rightarrow \left(
{{C(\lam)^t_0}}/{{\N\left(\lam\right)_0^t}}
\right)^*\rightarrow
{\left( C(\lam)^t_0\right)}^*
\rightarrow (\N\left(\lam\right)_0^t)^*
\rightarrow 0,
\end{align}
where ${}^*$ is defined in \eqref{eq:graded-dual}.
\begin{Pro}\label{Pro:Null-is-Null}
Let $\lam\in \dual{\h}$ be as in
Proposition \ref{Pro:later}.
Then,
\begin{align*}
\Wcoho^{i}\left(\n,M(\lam)^*\right)_a=
\Hom_{\C}
\left(H_{i}
\left(
{{C(\lam)^t_0}}/{{\N\left(\lam\right)_0^t}}
\right)_{a},\C\right)
\end{align*}
for all $i$ and $a\in \C$.
\end{Pro}
\begin{proof}
We first claim that
\begin{align}\label{eq:vanishing-null-dual}
H^{\bullet}
\left((\N\left(\lam\right)_0^t)^*
\right)\equiv 0.
\end{align}
Considering the spectral sequence
described in \cite[3.2]{FKW},
it is enough to show
that
$H^{\bullet}
\left((\N\left(\lam\right)_0^t)^*,
d_{\n}^{\st}\right)\equiv 0$.
But since
the action of
$d_{\n}^{\st}$ is compatible
with the weight space decomposition,
this is equivalent
to  $H_{\bullet}
\left({{\N\left(\lam\right)_0^t}},d_{\n^t}^{\st}
\right)\equiv 0$.
Hence Proposition \ref{Pro:later}
proves \eqref{eq:vanishing-null-dual}.

Now consider the
long exact sequence induced by
\eqref{eq:exact}.
Then,
by
\eqref{eq:vanishing-null-dual},
we get
$H^{\bullet}\left(\left(C(\lam)^t\right)^*\right)
=H^{\bullet}
\left(\left(
{{C(\lam)^t_0}}/{{\N\left(\lam\right)_0^t}}
\right)^*\right)$.
But Proposition \ref{Pro:finite-dim-of-C/N}
implies
\begin{align*}
H^{i}
\left(\left(
{{C(\lam)^t_0}}/
{{\N\left(\lam\right)_0^t}}
\right)^*\right)_a
=\Hom_{\C}\left(H_i\left(
{{C(\lam)^t_0}}/
{{\N\left(\lam\right)_0^t}}
\right)_{a}
,\C\right),
\end{align*}
for $a\in \C$.
Thus,
Proposition \ref{Pro:reduction-of-dual}
proves the proposition.
\end{proof}
\subsection{The cohomology
$\Wcoho^{\bullet}\left(\n,M(\lam)^*\right)$}
\begin{Pro}\label{Pro:iso-C/N}
Let $\lam\in \dual{\h}$.
\begin{enumerate}
\item
Let
$\n=\Ln{+}$.
For all $i\in \Z$ and $a\in \C$,
\begin{align*}
H_i
\left(
{{C(\lam)^t_0}}/{\N(\lam)_0^t}
\right)_a
\cong
\Wcoho^{-i}\left(\Ln{-},M(t_{-\srho\che}\circ
\lam)\right)_a.
\end{align*}
\item
Let
$\n=\Ln{-}$.
For all $i\in \Z$ and $a\in \C$,
\begin{align*}
H_i
\left(
{{C(\lam)^t_0}}/{\N(\lam)_0^t}
\right)_a
\cong
\Wcoho^{-i}\left(\Ln{-},M(w_0\circ
\lam)\right)_a.
\end{align*}
\end{enumerate}
\end{Pro}
\begin{proof}
(1)
By %Proposition \ref{Pro:iso-map-decomposition},
Proposition \ref{Pro:dec-of-complex-Verma}
and Proposition \ref{Pro:decompsotion-FB} (2),
we have
\begin{align*}
\Wcoho^{\bullet}\left(\n^t,M(t_{-\srho\che}\circ
\lam)\right)=
H^{\bullet}(B\*_{B_+}\C |t_{-\srho\che}\circ
\lam\ket,d_{\n^t} ).
\end{align*}
Observe
that
\begin{align}\label{eq:diffrential-shift1}
d_{\n}^t\phi_{t_{\srho\che}}(v)=
\phi_{t_{\srho\che}}(d_{\n^t}v),\quad
\DW_{\n}\phi_{t_{\srho\che}}(v)=
\phi_{t_{\srho\che}}(\DW_{\n^t}v)
\end{align}
for $v\in
C(\n^t,M(t_{-\srho\che}\circ\lam))$,
see Proposition \ref{Pro:tristing-heisenbeg}.
Therefore,
\begin{align*}
\Wcoho^{\bullet}\left(\n^t,M(t_{-\srho\che}\circ
\lam)\right)_a\cong
H^{\bullet}\left(
\phi_{t_{\srho\che}}
\left(B\*_{B_+}\C|t_{-\srho\che}\circ
\lam\ket \right),
d_{\n}^t\right)_a.
\end{align*}
for $a\in \C$.
Moreover,
the action of
$\widehat{J}_{\alpha}(n)$,
($\alpha\in \sproots$,
$n\in \Z$)
is well-defined on
$\phi_{t_{\srho\che}}\left(
C(\n^t,M(t_{-\srho\che}\circ\lam))\right)$,
and
we have
\begin{align}
\widehat{J}_{\alpha}(n)\phi_{t_{\srho\che}}(v)=
\phi_{t_{\srho\che}}(\widehat{J}_{\alpha}(n
+\height \alpha)v)\quad
(\alpha\in \sroots, n\in \Z).
\label{twisting_1}
\end{align}
Thus,
by Proposition
\ref{Pro:tristing-heisenbeg}
and \eqref{twisting_1}, it follows that
\begin{align*}
\phi_{t_{\srho\che}}\left(B\*_{B_+}\C
|t_{-\srho\che}\circ
\lam\ket\right)
=B\*_{t_{\srho\che}(B_+)}
\C\overline{|\lam\ket}.
\end{align*}
Here,
$\C\overline{|\lam\ket}$
is the
one-dimensional representation
of $t_{\srho\che}(B_+)$
appeared in \eqref{eq:C/N-+}.
Then,
by \eqref{eq:one-dim-+} and \eqref{eq:C/N-+},
we conclude
$
C(\lam)^t_0/\N(\lam)_0^t
=t_{\srho\che}\left(B\*_{B_+}\C |t_{-\srho\che}\circ
\lam\ket\right)
$
as complexes.
(2) can be similarly proved using
\begin{align}\label{eq:diffrential-shift2}
&d_{\n}^t\phi_{w_0}(v)=
\phi_{w_0}(d_{\n}v),\quad
\DW_{\n^t}\phi_{w_0}(v)=
\phi_{w_0}(\DW_{\n}v),\\
&\widehat{J}_{\alpha}(n)\phi_{w_0}(v)=
\phi_{w_0}(\widehat{w_0(J_{\alpha})}(n
)v)\quad
(\alpha\in \sroots, n\in \Z).
\label{twisting_2}
\end{align}
for $v\in
C(\n,M(w_0\circ\lam))$.
\end{proof}

\begin{Th}\label{Th:vanishin-dual}
Let $\lam\in \dual{\h}$.

{\rm (1)}
Suppose
$\bra \lam+\rho,\alpha\che\ket\not\in \Z_{\geq 1}$
for all
$\alpha\in \prroots\cap t_{\srho\che}(\nrroots)$.
Then, for all $a\in \C$,
\begin{align*}
\Wcoho^{i}\left(\Ln{+},M(\lam)^*\right)_a
\cong
\begin{cases}
\Hom_{\C}
\left(
\Wcoho^{0}\left(\Ln{-},M(t_{-\srho\che}\circ\lam)\right)_a,
\C\right)&(i=0)\\
\{0\}&(i\ne 0).
\end{cases}
\end{align*}

{\rm (2)}
Suppose $\bra \lam+\rho,\alpha\che\ket\not\in
\Z_{\geq 1}$
for all $\alpha\in \sproots$.
Then,
for all $a\in \C$,
\begin{align*}
\Wcoho^{i}\left(\Ln{-},M(\lam)^*\right)_a
\cong
\begin{cases}
\Hom_{\C}
\left(
\Wcoho^{0}\left(\Ln{-},M(w_0\circ\lam)\right)_a,
\C\right)&(i=0)\\
\{0\}&(i\ne 0).
\end{cases}
\end{align*}
\end{Th}
\begin{proof}
follows from
Proposition \ref{Pro:Null-is-Null},
Proposition
\ref{Pro:iso-C/N} and
Theorem \ref{Th:vanishing.for.Verma}.
\end{proof}
\subsection{Proof of
Proposition
\ref{Pro:later} (1)}
\label{subsection:The spectral seqeunce I}$ $

{\bf Step 1}
\quad Define the subspace
$\bar{F}^{p}{\bar{C}(\lam)^t}$,
$p\leq 0$,
of ${\bar{C}(\lam)^t}$ by
\begin{align*}
\bar{F}^{-p}{\bar{C}(\lam)^t}=
\bigoplus_{a\geq  \Whw{\n}{\lam}+p}
{\bar{C}(\lam)^t_{a}}
\subset {\bar{C}(\lam)^t}.
\end{align*}
Then,
by \eqref{eq:d-weights-of-bar-C-+},
\begin{align*}
&\dots \subset
\bar{F}^{-p} {\bar{C}(\lam)^t}
%\subset
%\bar{F}^{-p+1} {\bar{C}(\lam)^t}
\subset \dots
\subset
\bar{F}^0 {\bar{C}(\lam)^t}=
{\bar{C}(\lam)^t},\quad \bigcap_p\bar{F}^p
{\bar{C}(\lam)^t}=\{0\}.
\end{align*}
Notice that
$P\left(\bar{F}^{p}{\bar{C}(\lam)^t}\right)
\subset
\{\mu\in \dual{\h}_{\kappa}
;
\bra
\mu-\lam,\srho\che\ket \geq
-p\}$.
Thus,
\begin{align}
&\bar{F}^{p}{\bar{C}(\lam)^t} \subset
{\bar{\N}(\lam)^t}
\quad
\text{ for $p\leq -1$},\\
&{\bar{\N}(\lam)^t}=\left({\bar{\N}(\lam)^t}
\cap \bar{F}^{-1} {\bar{C}(\lam)^t}\right)
\+\sum\limits_{\bra \lam-\mu,\srho\che\ket<
0} \left({\bar{C}(\lam)^t}\right)^{\mu}_{\Whw{\n}{\lam}},
\label{eq:bar-N-dec-+}
\end{align}
where
$\left({\bar{C}(\lam)^t}\right)^{\mu}_{a}=
\left({\bar{C}(\lam)^t}\right)^{\mu}
\cap {\bar{C}(\lam)^t_{a}}$.

Define the subspace
$F^{p}{{C(\lam)^t_0}}
$,
$p\leq 0$,
of ${{C(\lam)^t_0}}$
by
\begin{align*}
F^p{{C(\lam)^t_0}}=B
\*_{t_{\srho\che}(B)
}\bar{F}^p{\bar{C}(\lam)^t}.
\end{align*}
Then,
\begin{align}
&\dots \subset
F^{-p} {{C(\lam)^t_0}}
\subset \dots
\subset
F^0 {{C(\lam)^t_0}}=
{{C(\lam)^t_0}},\quad
\bigcap_pF^p {{C(\lam)^t_0}}=\{0\},
\nonumber \\
&P\left({F}^{p}{C(\lam)^t}\right)\subset
\{\mu\in \dual{\h}_{\kappa}
;
\bra
\mu-\lam,\srho\che\ket \geq
-p\}
\label{eq:weight-of-Fdash},\\
& {F}^{p}{C(\lam)^t}
\subset {{\N\left(\lam\right)_0^t}}
\quad \text{for $p\leq -1$}.
\label{eq:F-p-is-sub-of-N}
\end{align}
Let
$F^p({{C(\lam)^t_0}})^{\mu}=
F^p{{C(\lam)^t_0}}\cap
({{C(\lam)^t_0}})^{\mu}$.
Then,
by \eqref{eq:weight-of-Fdash},
for a given
$\mu\in \dual{\h}$,
$F^p({{C(\lam)^t_0}})^{\mu}=\{0\}$
for $p\ll 0$.
Hence,
$\{{F}^{p}({C(\lam)^t})^{\mu}\}$
defines a convergent filtration
bounded below on each subcomplex
$({C(\lam)^t})^{\mu}$,
$\mu\in \dual{\h}$.
\begin{Pro}\label{Pro:first_filtration}
Let $\lam\in \dual{\h}$.
\begin{enumerate}
\item
$d_{\n^t}^{\st}\bar{F}^p{\bar{C}(\lam)^t}\subset
\bar{F}^p{\bar{C}(\lam)^t}
+F^{p-1} {{C(\lam)^t_0}}$.
\item
$
d_{\n^t}^{\st}F^p {{C(\lam)^t_0}}
\subset F^p {{C(\lam)^t_0}}$.
\end{enumerate}
\end{Pro}
\begin{proof}
(1) follows from
the commutativity
of $d_{\n^t}^{\st}$
with $\DW_{\n}$
and
the fact that
the operators
$\widehat{J}_{\alpha}(-n)$,
$\psi_{\alpha}(-n)$
(
$\alpha\in \sproots$,
$n>\height \alpha$)
and $\widehat{h}(-n)$
($h\in \sh$,
$n>0$)
have negative eigenvalues
with respect to the
adjoint
action of $\DW_{\n}$.
(2) follows
from
(1) and the definition of
$F^p {{C(\lam)^t_0}}$.
\end{proof}

Consider the spectral sequence
${E}^r\Rightarrow H_{\bullet}
\left({{C(\lam)^t_0}},d_{\n^t}^{\st}\right)$
corresponding
to the filtration
$\{F^p{{C(\lam)^t_0}}\}$.
We have:
${E}_{p,\bullet}^1=
H_{\bullet}\left(F^p{{C(\lam)^t_0}}
/F^{p-1}{{C(\lam)^t_0}},
d_{\n^t}^{\st}
\right)$.

Let
\begin{align*}
&{C}(\lam)'=
\sum_p
F^p
{{C(\lam)^t_0}}/F^{p-1}{{C(\lam)^t_0}},\\
&
\N(\lam)'=\im:
\sum
F^p
{{\N\left(\lam\right)_0^t}}/
F^{p-1}
{{\N\left(\lam\right)_0^t}}
\hookrightarrow C(\lam)',
\end{align*}
where
$F^p
{{\N\left(\lam\right)_0^t}}=
{{\N\left(\lam\right)_0^t}}\cap F^p
{{C(\lam)^t_0}}$.

\begin{Pro}\label{Pro:intermidiate}
$d_{\n^t}^{\st}\N(\lam)'
\subset \N(\lam)'$,
and if $\lam\in \dual{\h}$
satisfies
$\bra
\lam+\rho,\alpha\che\ket\not\in
\Z_{\geq 1}$
for all $\alpha\in \prroots\cap t_{\srho\che}
(\nrroots)$,
then
$H_{\bullet}
\left(\N(\lam)',d_{\n^t}^{\st}
\right)\equiv 0$.
\end{Pro}
Proposition
\ref{Pro:intermidiate}
will be proven in Step 3.
Note that Proposition
\ref{Pro:intermidiate}
implies
Proposition
\ref{Pro:later} (1).
Indeed,
by
\eqref{eq:bar-N-dec-+},
\eqref{eq:F-p-is-sub-of-N}
and Proposition \ref{Pro:first_filtration},
$d_{\n^t}^{\st}\N(\lam)'
\subset \N(\lam)'$
implies
$d_{\n^t}^{\st}{{\N\left(\lam\right)_0^t}}
\subset {{\N\left(\lam\right)_0^t}}$,
and
by
\eqref{eq:F-p-is-sub-of-N} again,
$H_{\bullet}
\left(\N(\lam)',d_{\n^t}^{\st}
\right)\equiv 0$
implies
that
${E}^r$
degenerates at the
$E_1$-term itself
and that
$
H_{\bullet}\left({{C(\lam)^t_0}},
d_{\n^t}^{\st}
\right)
=H_{\bullet}\left({C(\lam)^t}
/\N(\lam)^t,
d_{\n^t}^{\st}
\right)
$,
that is,
$H_{\bullet}\left(\N(\lam)^t,
d_{\n^t}^{\st}
\right)\equiv 0$.

\medskip

{\bf Step 2}
\quad
Let
\begin{align*}
&
\bar{C}(\lam)'=\im:
\sum_p \bar{F}^p{\bar{C}(\lam)^t}/
\bar{F}^{p-1}{\bar{C}(\lam)^t}
\hookrightarrow
{C}(\lam)'.
\end{align*}
Then,
$C(\lam)'=B\*_{t_{\srho\che}(B)
}{\bar{C}(\lam)'}$
and
$\bar{C}(\lam)'$ is a subcomplex of
${C}(\lam)'$
by Proposition \ref{Pro:first_filtration} (1).
Observe  that by  definition,
it is the following quotient complex of
${C(\lam)^t}$:
\begin{align}\label{bar-C-is-quotient}
\bar{C}(\lam)'
={C(\lam)^t}/
\haru\left\{\begin{array}
{l}
\widehat{J}_{\alpha}(-n)v,
\psi_{\alpha}(-n)v,
\widehat{J_i}(-m)v,\\
;\alpha\in \sproots,
n>\height \alpha,
i\in \sI,
m>0,v\in {C(\lam)^t}
\end{array}\right\}.
\end{align}

\begin{Pro}\label{Pro:coho-C-bar}$ $
$H_i\left(\bar{C}(\lam)',
d_{\n^t}^{\st}\right)=
H^{\semiinf-i}\left(t_{\srho\che}(\g_-),
M(\lam)\right)$.
\end{Pro}
\begin{proof}
By the duality of
the standard semi-infinite cohomology (\cite{Feigin}),
we have
\begin{align}\label{eq:duality-of-Feigin}
H^{\semiinf-i}\left(t_{\srho\che}(\g_-),
M(\lam)\right)
=\coho{i}
\left(t_{\srho\che}(\g_+),M(\lam)^*
\right)^*.
\end{align}
Thus,
it is sufficient
to show that
$H^{\bullet}
\left((\bar{C}(\lam)')^*
\right)
=\coho{\bullet}
\left(t_{\srho\che}(\g_+),M(\lam)^*
\right)$.

Let
$\m=
t_{\srho\che}(\g_+)\cap \g_-
=t_{\srho\che}(\n_+)\cap \n_-$.
Since
$M(\lam)^*$
is $t_{\srho\che}(\g_+)\cap \g_+$-cofree,
by \cite[Theorem 2.1]{V1}
it follows that
\begin{align*}
&\coho{i}
\left(t_{\srho\che}(\g_+),M(\lam)^*
\right)\\
&=
\begin{cases}
\{0\}&(i>0),\\
H_{-i}\left(
\Big(M(\lam)^*\* \Lam
(\m)\Big)^{t_{\srho\che}(\g_+)\cap \g_+}, d_{\m}^{\st}
\right)
&(i\leq 0)
\end{cases}
\end{align*}
where
$t_{\srho\che}(\g_+)\cap \g_+$
acts on $\Lam (\m)$
via the identification
$\m=t_{\srho\che}(\g_+)/
\left(t_{\srho\che}(\g_+)\cap \g_+\right)$
and
$d_{\m}^{\st}$ is the differential of $\m$-homology,
i.e,
\begin{align*}
d_{\m}^{\st}=
\sum_\ud{\alpha\in \sproots,
}{ -\height \alpha\leq n<0}
J_{\alpha}
(n)
\psi_{-\alpha}(-n)
-\frac{1}{2}\sum_\ud{\alpha,\beta,\gamma\in
\sproots
}{
-\height \alpha\leq k<0,
-\height \beta \leq l<0}c_{\alpha,\beta}^{\gamma}
\psi_{-\alpha}(-k)\psi_{-\beta}(-l)\psi_{\gamma}
(k+l).
\end{align*}

On the other hand,
by \eqref{bar-C-is-quotient}
and \eqref{eq:transpose-hat-J},
we have
\begin{align}\label{eq:identificaion-of-D-bar-C}
(\bar{C}(\lam)')^*
=
\Big(M(\lam)^*\*
\Lam (\m)\Big)^{ \t_{\srho\che}(\g_+)
\cap \g_+}\subset
\Big(M(\lam)^*\*
\Lam (\n_-)\Big)^{ \n_+}.
%=C(\lam)^0.
\end{align}
This can be proved in the same way as
Proposition \ref{Pro:identification-of-C-top-0}.
Therefore,
by Proposition
\ref{Pro:identification-of-C-top-0},
it is now sufficient
to check that $d_{\n_-}^{\st}$ acts
as $d_{\m}^{\st}$ on
the right-hand-side of
\eqref{eq:identificaion-of-D-bar-C}.
But this is easy to see.
\end{proof}
{\bf Step 3}~~
Define the subspace
$\Fdash^pC(\lam)'$,
$p\leq 0$,
of $C(\lam)'$
by
\begin{align*}
\Fdash^pC(\lam)'=B\*_{t_{\srho\che}(B_+)
}\bar{\Fdash}^p{\bar{C}(\lam)'},
\end{align*}
where
\begin{align*}
&\bar{\Fdash}^p{\bar{C}(\lam)'}=
\bigoplus_\ud{\mu}{
\bra \mu-\lam,\srho\che\ket\geq -p}
\bar{\Fdash}^p{(\bar{C}(\lam)')}^{\mu}
\subset
{\bar{C}(\lam)'}.
\end{align*}
Then,
similarly as in the step 1,
we have:
\begin{align}
&\N(\lam)'
=
\Fdash^{-1}C(\lam)'
\label{eq:N-and-F-1-+}
\\
&\dots
\subset \Fdash^{-p}C(\lam)'
\subset \dots
\subset
\Fdash^0C(\lam)'=
C(\lam)',
\quad \bigcap \Fdash^p
C(\lam)'=\{0\},\\
&\Fdash^{p}(C(\lam)')^{\mu}=\{0\}
\quad (p\ll 0).
\label{eq:fdd-sub}
\end{align}
Since $d_{\n^t}^{\st}
\bar{\Fdash}^p
\bar{C}(\lam)'\subset\bar{\Fdash}^p
\bar{C}(\lam)'
$
and $t_{\srho\che}(B_+)\bar{\Fdash}^p
\bar{C}(\lam)'\subset\bar{\Fdash}^p
\bar{C}(\lam)'$,
it follows that
$ d_{\n^t}^{\st}\Fdash^{p}C(\lam)'
\subset \Fdash^{p}C(\lam)'$.
In particular,
$d_{\n^t}^{\st}\N(\lam)'
\subset \N(\lam)'$.

Let
${E'}^r\Rightarrow
H_{\bullet}
(C(\lam)',d_{\n^t}^{\st})$
be the corresponding
spectral sequence.
We have:
${E'}_{p,\bullet}^1=
H_{\bullet}\left(\Fdash^pC(\lam)'
/\Fdash^{p-1}C(\lam)',
d_{\n^t}^{\st}
\right)$.
Notice that
$\bar{\Fdash}^p
\bar{C}(\lam)'/\bar{\Fdash}^{p-1}
\bar{C}(\lam)'$
is a subcomplex
of
$\Fdash^pC(\lam)'
/\Fdash^{p-1}C(\lam)'$
and that
\begin{align*}
\bar{\Fdash}^p
\bar{C}(\lam)'/\bar{\Fdash}^{p-1}
\bar{C}(\lam)'=
\bigoplus\limits_\ud{\mu}{
\bra \mu-\lam,\srho\che\ket
=-p}(\bar{C}(\lam)')^{\mu}
\end{align*}
as a complex.

Consider
$\phi_{t_{\srho\che}}
\left(
B\*_{B_+}
\C |\mu\ket \right)
\subset \phi_{t_{\srho\che}}
\left(C(\n,M(\mu)
\right)$
as a complex with differential
$d_{\n}^t$
as in the proof of Proposition
\ref{Pro:iso-C/N}.
\begin{Pro}\label{Pro:complex-str-of-Fp-Fp-+}
Let $\lam\in \dual{\h}$.
\begin{align*}
\Fdash^pC(\lam)'
/\Fdash^{p-1}C(\lam)'=
\bigoplus_\ud{\mu}{
\bra\mu-\lam,\srho\che\ket=-p}
\phi_{t_{\srho\che}}
\left(
B\*_{B_+}
\C |t_{-\srho\che}
\circ \mu\ket \right)\*
(\bar{C}(\lam)')^{\mu}
\end{align*}
as a complex.
\end{Pro}
\begin{proof}
is similar to that of Proposition \ref{Pro:iso-C/N}.
Indeed,
we have
\begin{align*}
\Fdash^pC(\lam)'
/\Fdash^{p-1}C(\lam)'&=
B\*_{t_{\srho\che}(B_+)}
\left(\bar{\Fdash}^p{\bar{C}(\lam)'}
/\bar{\Fdash}^{p-1}{\bar{C}(\lam)'}\right)\\
&=\bigoplus_{\mu}
B\*_{t_{\srho\che}(B_+)}
\left(\bar{\Fdash}^p{\bar{C}(\lam)'}
/\bar{\Fdash}^{p-1}{\bar{C}(\lam)'}\right)^{\mu},
\end{align*}
and
$\left(\bar{\Fdash}^p{\bar{C}(\lam)'}
/\bar{\Fdash}^{p-1}{\bar{C}(\lam)'}\right)^{\mu}$
is a direct sum of
copies of $\C \overline{|\mu\ket}$
as a
$t_{\srho\che}(B_+)$-module.
\end{proof}
\begin{Pro}\label{Pro:one-more-step-to-go-+}
Let $\lam$ as in Proposition
\ref{Pro:later} (1).
Then,
\begin{align*}
\text{$H_{\bullet}
\left(\Fdash^pC(\lam)'
/\Fdash^{p-1}C(\lam)',d_{\n^t}^{\st}\right)
=\{0\}\quad $
$(p\ne 0)$.}
\end{align*}
\end{Pro}
\begin{proof}
By Theorem \ref{Th:cofreeness_of_Verma},
the assumption on  the weight $\lam$
implies
$M(\lam)$ is
is cofree over $t_{\srho\che}(\g_-)\cap
\g_+$.
Since $M(\lam)$
is obviously free over $t_{\srho\che}(\g_-)\cap
\g_-$,
\cite[Theorem 2.1]{V1}
implies
\begin{align*}
H^{\semiinf+i}\left(t_{\srho\che}
(\g_-),M(\lam)\right)^{\mu}=
\begin{cases}
\C_{\lam}&\text{($i=0$ and $\mu=\lam$)}\\
\{0\}&\text{(otherwise)}.
\end{cases}
\end{align*}
Thus Proposition \ref{Pro:coho-C-bar}
and Proposition
\ref{Pro:complex-str-of-Fp-Fp-+}
prove the proposition.
\end{proof}

By Proposition \ref{Pro:one-more-step-to-go-+},
${E'}^r$
degenerates at the
$E_1$-term itself,
i.e,
\begin{align*}
H_{\bullet}\left(C(\lam)',
d_{\n^t}^{\st}
\right)=H_{\bullet}\left({\Fdash}^0C(\lam)'
/{\Fdash}^{-1}C(\lam)',
d_{\n^t}^{\st}
\right)
=H_{\bullet}\left(C(\lam)'
/\N(\lam)',
d_{\n^t}^{\st}
\right).
\end{align*}
Here the last equality follows from \eqref{eq:fdd-sub}.
This implies
$H_{\bullet}\left(\N(\lam)',
d_{\n^t}^{\st}
\right)\equiv 0$.
Proposition
{\ref{Pro:intermidiate}}
is proved.
Thus,
Proposition
\ref{Pro:later} (1) is proved.
\renewcommand{\qedsymbol}{$\blacksquare$}
\qed
\renewcommand{\qedsymbol}{$\square$}

\subsection{Proof of
Proposition
\ref{Pro:later} (2)}
\label{subsection:spectral-I}
We omit the most of the
proof of (2).
Indeed,
its proof is simpler than
(1):
By \eqref{eq:d-weights-of-bar-C--},
step 1
in the previous section
is not needed for this case
and
the argument in step 2 is replaced
by the following proposition.
\begin{Pro}\label{Pro:dCbar-subbarC--}
Let $\lam\in \dual{\h}$.
\begin{enumerate}
\item
${\bar{C}(\lam)^t}=
\bar{M}(\bar{\lam})\* \Lam (\snp)^*
\subset M(\lam)\*\F
$.
Here,
$\bar{M}(\bar{\lam})$ is the Verma module
of $\sg$
of highest weight $\bar{\lam}$
identified with
$U(\snn)v_{\lam}
\subset M(\lam)$.
\item
$d_{\n^t}^{\st}{\bar{C}(\lam)^t}\subset
{\bar{C}(\lam)^t}$
and $$H_{i}
\left({\bar{C}(\lam)^t},d_{\n^t}^{\st}\right)
=H^{-i}\left(\snp,\bar{M}(\bar{\lam})\right),$$
where
${\bar{C}(\lam)^t}
=\sum_i
{\bar{C}^{i}(\lam)^t}$,
${\bar{C}^{i}(\lam)^t}=
{\bar{C}(\lam)^t}\cap
{C(\lam)^t}
$.

\end{enumerate}\end{Pro}
\begin{proof}
(1) follows from the fact that
$\widehat{J}_{-\alpha}(0)$
acts as ${J}_{-\alpha}(0)$
on
$M(\lam)\* \C\1$.
(2)
easily follows form (1).
Indeed,
\begin{align*}%\label{eq:diffrential-of-bar}
{d_{\n^t}^{\st}}_{|{\bar{C}(\lam)^t}}=
\sum_{\alpha\in \sproots}
J_{\alpha}(0){\psi}_{-\alpha}(0)
-\frac{1}{2}\sum_{\alpha,\beta,\gamma\in \sproots}
c_{\alpha,\beta}^{\gamma}{\psi}_{-\alpha}(0)
{\psi}_{-\beta}(0)\psi_{\gamma}(0).
\end{align*}
Thus,
by \cite[Proposition 4.7]{DGK},
it follows that
$H_{i}
\left({\bar{C}(\lam)^t},d_{\n^t}^{\st}\right)
=H^{-i}\left(\snp,\bar{M}(\bar{\lam})\right)$.
\end{proof}
\renewcommand{\qedsymbol}{$\blacksquare$}

\qed

\renewcommand{\qedsymbol}{$\square$}

\section{Estimate on $\DW_{\n}$-eigenvalues}
\label{section:estimate}
In this section
we shall give an
estimate
of $\DW_{\n}$-eigenvalues
of $\Wcoho^{\bullet}(\n,V)$
for $V\in \BGG_{\kappa}^{[ \Lam]}$
under the restriction
on $\Lam$
as in Introduction.
The results in this section
will be needed  when
$\kappa\in \Q_{>0}$.
Let
$\fd_{\kappa,+}
=\{ \Lam\in \dual{\h}_{\kappa};
\bra \Lam+\rho,\alpha\che\ket
\geq 0
\text{ for all }\alpha\in R^{\Lam}_+\}$,
the set of dominant weights of level $\kappa-h\che$.
Then,
$
\BGG_{\kappa}=\bigoplus\limits_{\Lam
\in \fd_{\kappa,+}}\BGG_{\kappa}^{[\Lam]}
$
if $\kappa\not\in \Q_{\leq 0}$.
\subsection{The use of the standard
semi-infinite cohomology
}
In
\cite[3.2]{FKW},
it was shown that
there exists a converging spectral sequence
$E_{r}^{p,q}\Rightarrow \Wcoho^{\bullet}(\n,V) $,
$V\in \BGG_{\kappa}$,
such that
$E^{\bullet,q}_1=\coho{q}(\n,V) $
and the corresponding filtration is
compatible with the action of
$\DW_{\n}$.
Let
$\coho{i}(\n,V)_{a}$
be the  eigenspace
of $\DW_{\n}$
with the eigenvalue $a\in \C$.
Then,
%by \eqref{eq:description_of-eigen-values},
\begin{align}
&\coho{i}(\n,V)_{a}=
\bigoplus_\ud{\lam\in
\dual{\h}}{
\Whw{\n}{\lam}=a}{\coho{i}(\n,V)}^{\lam}.
\label{eq:DW-eigenvalues-of-standard}
\end{align}
The following proposition
is clear.
\begin{Pro}\label{Pro:gr}
$ \Wcoho^{i}\left(\n,V\right)_a\ne \{0\}$
only if
$\coho{i}(\n,V)_{a}\ne \{0\}$
{\rm{(}}$V\in \BGG_{\kappa}$,
$i\in \Z$,
$a\in \C$ {\rm{)}}.
\end{Pro}

\subsection{The
formal character}
For $\lam\in \dual{\h}$,
let $I(\lam)$ be the irreducible
representation of $\H$
of highest weight $\lam +h\che\Lam_0$.
Since the category of highest weight $\H$-modules is completely
reducible,
$C(\n,V)$, $V\in \BGG_{\kappa}$,
decomposes into a direct sum of
$I(\lam)$:
$
C(\n,V)
\cong
\bigoplus_{\lam\in \dual{\h}_{\kappa}
} \B_{\lam}^{\bullet}(\n,V)
\* \Hsim{\lam}
$.
Here,
\begin{align}
&\B_{\lam}^{i}(\n,V)\nonumber\\&=
\{v\in C^{i}(\n,V);
\widehat{h}(n)\cdot v=0,~
\widehat{h}(0)\cdot v=\lam(h)v,~
(h\in \sh,n>0),~
\Dtot\cdot v=\lam(\Dg)v\}
,\nonumber\\
&=\{v\in C^{i}(\n,V)^{\lam};\widehat{h}(n)\cdot v=0~
(h\in \sh,n>0)\}.\label{eq:branching}
\end{align}
Note that $\dim \B_{\lam}^{\bullet}(\n,V)<\infty$
by definition.
By the commutativity
of $d_{\n}^{\st}$ with  the action of
$\H$,
it follows that
\begin{align}\label{eq:branching_of_cohomology}
\coho{\bullet}(\n,V)
\cong
\bigoplus_{\lam\in \dual{\h}_{\kappa}
}
H^{\bullet}(\B_{\lam}^{\bullet}(\n,V))
\* \Hsim{\lam}\quad
(V\in \BGG_{\kappa}).
\end{align}
Here,
$H^{\bullet}(\B_{\lam}^{\bullet}(\n,V))
=H^{\bullet}(\B_{\lam}^{\bullet}(\n,V),d_{\n}^{\st})$.
\begin{Rem}
By \cite{HT},
it follows that
the sum in the
right-hand-side in \eqref{eq:branching_of_cohomology}
is taken over
$\lam\in \dual{\h}_{\kappa}$ such that
$|\lam+\rho|^2=|\Lam+\rho|^2$
for $V\in \BGG_{\kappa}^{[\Lam]}$.
\end{Rem}
The following is clear by
\eqref{eq:DW-eigenvalues-of-standard},
\eqref{eq:branching_of_cohomology}
and
Proposition \ref{Pro:gr}.
\begin{Lem}\label{Lemma:for-estimate}
The $\DW_{\n}$-eigenvalues
of
$\Wcoho^{i}\left(\n,V\right)$,
$V\in
\BGG_{\kappa}$,
are contained in the set
\begin{align*}%\label{eq:weights-in-contained-first}
\bigcup_\ud{\lam\in \dual{\h}
}{ H^{i}(\B_{\lam}^{\bullet}(\n,V))\ne
\{0\}}\Whw{\n}{\lam}-\Z_{\geq 0}.
\end{align*}
\end{Lem}

Define
$\ch  H^{\bullet}(\B_{\lam}^{\bullet}(\n,V))
=\sum_{i\in \Z} z^i
\dim H^{i}(\B_{\lam}^{\bullet}(\n,V))
$.
Then,
by \eqref{eq:branching_of_cohomology},
\begin{align}\label{eq:ch-B-and-Ch-C}
\sum_{\lam}\ch H^{\bullet}(\B_{\lam}^{\bullet}(\n,V))
e^{\lam}=\prod_{\alpha\in \piroots}
(1-e^{-\alpha})^{\dim \g_{\alpha}}
\ch \coho{\bullet}(\n,V),
\end{align}
where
$\ch \coho{\bullet}(\n,V)
=\sum_{i\in \Z}z^i \sum_{\lam\in
\dual{\h}}e^{\lam}
\dim \coho{i}(\n,V)^{\lam}$.
%Similarly,
%set
%$\ch C(\n,V)=
%\sum_{i\in \Z}z^i \sum_{\lam\in
%\dual{\h}}e^{\lam}
%\dim C(\n,V)^{\lam}$.
\subsection{The estimate on $\DW_{\n}$-eigenvalues}
\begin{Lem}\label{Lem:lem-on-roots}
Let $\Lam\in \dual{\h}$.
For $w\in \eW$,
the following conditions are equivalent
\begin{enumerate}
\item
$\bra \Lam+\rho,\alpha\che\ket
\not\in \Z$
for all $\alpha\in \prroots\cap w(\nrroots)$.
\item $w\inv (R^{\Lam}_+)\subset \prroots$.
\end{enumerate}

\end{Lem}
\begin{proof}
(1) is equivalent to
$R^{\Lam}_+\cap \prroots\cap w(\nrroots)
=\emptyset$.
On the other hand,
(2) is equivalent to $R^{\Lam}_+\subset
\prroots\cap w(\prroots)$.
But these two conditions are
equivalent.
\end{proof}

\begin{Lem}\label{Lem:estimate-height}
Let $\Lam\in \dual{\h}_{\kappa}$ such that
\begin{align*}
&\begin{cases}
\text{$\bra \Lam+\rho,\alpha\che\ket\not\in \Z$
for all
$\alpha\in \prroots\cap t_{\srho\che}(\nrroots)$.}
&(\text{if $\n=\Ln{+})$}
\\
\text{$\bra \Lam+\rho,\alpha\che\ket\not\in \Z$
for all
$\alpha\in \sproots$.}
&(\text{if $\n=\Ln{-})$}.
\end{cases}
\end{align*}
Then,
$\Whw{\n}{\Lam}-\Whw{\n}{\mu}
\in  \Z_{\geq \height_{\Lam}
(\Lam-\mu)}$
for all $\mu\in \Lam-Q_+^{\Lam}$.
\end{Lem}
\begin{proof}
Let $\n=\Ln{+}$.
Notice that
$t_{-\srho\che}(\prroots)
\cap \sproots=\emptyset$.
Thus,
by the assumption
and Lemma \ref{Lem:lem-on-roots},
$t_{-\srho\che}(\alpha)\in
\prroots\backslash \sproots$
for all $\alpha\in \Pi^{\Lam}$.
Hence,
$\bra t_{-\srho\che}(\alpha),\Dg\ket\in \Z_{\geq 1}$
for all $\alpha\in \Pi^{\Lam}$,
and thus,
$\bra t_{-\srho\che}(\mu),\Dg\ket
\in \Z_{\geq \height_{\Lam}(\mu)}$
for $\mu \in Q^{\Lam}_+$.
But
\begin{align}\label{eq:diffrence-of-weights}
\Whw{\n}{\Lam}-\Whw{\n}{\mu}
=\bra \Lam-\mu,\srho\che+\Dg\ket
=\bra t_{-\srho\che}( \Lam-\mu),\Dg\ket.
\end{align}
Therefore the assertion follows.
The $\n=\Ln{-}$   case  follows from the formula
\begin{align}\label{eq:diffrence-of-weights-}
\Whw{\n}{\Lam}-\Whw{\n}{\mu}
=\bra \Lam-\mu,\Dg\ket.
\end{align}
and the fact that
$R^{\Lam}_+\subset \prroots\backslash \sproots$.
\end{proof}

\begin{Pro}\label{Pro:pre-crutial-estimate}
Let $\Lam\in \dual{\h}_{\kappa}$
be as in Lemma \ref{Lem:estimate-height}.
Then,
for all
$i\in \Z$
and $V\in \BGG_{\kappa}^{[\Lam]}$,
$\DW_{\n}$-eigenvalues of
$\Wcoho^{i}\left(\n,V\right)$
is contained in the set
\begin{align*}
\bigcup_\ud{\mu\in
W^{\Lam}\circ \Lam
}{
[V:L(\mu)]\ne 0} \Whw{\n}{\mu}-\Z_{\geq
|i|}.
\end{align*}
Here,
$[V:L(\mu)]$
is the multiplicity
of $L(\mu)$ in $V$
in the sense of \cite{DGK}.
\end{Pro}
Proof of Proposition \ref{Pro:pre-crutial-estimate}
is given at the end of this section.

Lemma \ref{Lem:estimate-height} and
Proposition \ref{Pro:pre-crutial-estimate}
imply:
\begin{Co}\label{CO:estimate-of-weights-kappa-Q-+}
Let $\kappa\in \Q_{> 0}$.
Suppose that $\Lam\in \fd_{\kappa,+}$
satisfies the condition in Lemma
\ref{Lem:estimate-height}. Then,
for all
$i\in \Z$
and $V\in \BGG_{\kappa}^{[\Lam]}$,
$\DW_{\n}$-eigenvalues of
$\Wcoho^{i}\left(\n,V\right)$
is contained in the set
$\Whw{\n}{\Lam}-\Z_{\geq
|i|}$.
\end{Co}
\begin{Lem}\label{Lem:res-from-hws}
Let $\Lam\in  \fd_{\kappa,+}$,
$\kappa\in \Q_{>0}$.
Then,
for a given
$N\in \Z_{\geq 0}$
and $V\in \BGG_{\kappa}^{[\Lam]}$,
there exists a finitely generated submodule
$M$ of $V$
such that
$[V/M:L(\mu)]=0$ if $\height_{\Lam}(\Lam-\mu)\leq
N$.
\end{Lem}
\begin{proof}
Let
$\{0\}=V_0\subset V_1\subset V_2\subset \dots $
be a highest weight series of $V$,
that is,
a filtration
of $V$ such
that
(1) $V=\bigcup V_i$,
(2) Each
subquotient
$V_i/V_{i-1}$ is a quotient of $M(\mu_i)$
for some $\mu_i\in \dual{\h}$,
and (3) $\mu_j-\mu_i\not\in Q_+$
for $i<j$.
Since $V\in \BGG_{\kappa}^{[\Lam]}$,
it follows that
$\mu_i\in W^{\Lam}\circ \Lam$
for all $i$.
We may assume that
$V_i\ne \{0\}$ for all $i$,
because there is nothing to show
if $V$ is finitely generated.
Since
$\{\lam\in W^{\Lam}\circ \Lam;
\height_{\Lam}(\Lam-\lam)\leq
N
\}$ is a finite set,
there exists an integer $k$
such that
$\height_{\Lam}(\Lam-\mu_i)>N$
for all $i>k$.
Let $M=V_k$.
Then,
$P(V/M)\subset \bigcup_{i>k}\mu_i-Q_+$,
and therefore,
$[V/M:L(\mu)]=0$
if $\height_{\Lam}(\Lam-\mu)\leq
N$.
\end{proof}
\begin{Pro}\label{Pro:crutial-estimate}
Let $\Lam\in  \fd_{\kappa,+}$,
$\kappa\in \Q_{>0}$.
Suppose that
$\Lam$ satisfies
satisfies the condition in Lemma
\ref{Lem:estimate-height}.
Let $V\in \BGG_{\kappa}^{[\Lam]}$
and suppose $a\in \C$ is given.
\begin{enumerate}
\item There exists a finitely generated
submodule $M$
of $V$ such that
$\Wcoho^{\bullet}(\n,V)_a\cong
\Wcoho^{\bullet}(\n,M)_a$.
\item
There exists a
quotient $M'$
of $V$ such that
$(M')^*$ is finitely generated and
$\Wcoho^{\bullet}(\n,V)_a\cong
\Wcoho^{\bullet}(\n,M')_a$.
\end{enumerate}
\end{Pro}
\begin{proof}
By Corollary \ref{CO:estimate-of-weights-kappa-Q-+},
we may assume that
$a\in  \Whw{\n}{\Lam}-\Z_{\geq 0}$.
Let $N=\Whw{\n}{\Lam}-a$.

(1)
By Lemma \ref{Lem:res-from-hws},
there exists a  finitely generated
submodule $M$ of $V$ such that
$[V/M:L(\mu)]=0$ if $\height_{\Lam}(\Lam-\mu)\leq
N$.
Then,
by Lemma \ref{Lem:estimate-height}
and Proposition
\ref{Pro:pre-crutial-estimate},
it follows that
\begin{align}\label{eq:temp-in-Pro}
\Wcoho^{\bullet}(\n,V/M)_{a}=\{0\}
\quad (a\geq  \Whw{\n}{\Lam}-N)
\end{align}
Consider the exact sequence
$0\rightarrow M\rightarrow V\rightarrow
V/M\rightarrow 0$ in $\BGG_{\kappa}^{[\Lam]}$.
It induces
the long exact sequence of semi-infinite
cohomology.
Clearly,
its restriction
to  a $\DW_{\n}$-eigenspace
remains exact.
Thus,
(1) follows from
\eqref{eq:temp-in-Pro}.
(2) is similarly proved as  (1).
Indeed,
let $M$ be a finitely generated submodule
of $V^*$ such that
$[V^*/M:L(\mu)]=0$ if $\height_{\Lam}(\Lam-\mu)\leq
N$.
Then,
$0\rightarrow (V^*/M)^*\rightarrow
V\rightarrow M^*\rightarrow 0$
and
$[(V^*/M)^*:L(\mu)]=0$ if $\height_{\Lam}(\Lam-\mu)\leq
N$.
\end{proof}

\subsection{Proof of Proposition
\ref{Pro:pre-crutial-estimate}}
Let
$\w_{\n}=
\begin{cases}
t_{\srho\che}&(\text{if $\n=\Ln{+})$}\\
w_0&(\text{if $\n=\Ln{-})$}
\end{cases}$
and $\m=\g_{\w_{\n}}^t$.
Thus,
$\m=\begin{cases}
t_{\srho\che}(\g_+)\cap
\g_-&(\text{if $\n=\Ln{+})$}\\
\sn_-&(\text{if $\n=\Ln{-})$}.
\end{cases}$
Note
$\m\subset \n_-$ for the either case.

Let $\Lam\in
\dual{\h}_{\kappa}$ be as in Lemma
\ref{Lem:estimate-height}. Then,
any objects in $\BGG_{\kappa}^{[\Lam]}$
is free over $\m$
by
Theorem \ref{Th:global_cofreeness}.
Therefore,
\begin{align}\label{eq:ch-of-m-coinvariants}
\ch (V/\m V)=
\prod_{\alpha\in \rroots(\m)}
(1-e^{\alpha})\ch V
\quad (V\in \BGG_{\kappa}^{[\Lam]}).
\end{align}
For
$V\in \BGG_{\kappa}^{[\Lam]}$,
$\Lam\in \dual{\h}_{\kappa}$,
define
$[V:M(\mu)]\in \Z$,
$\mu\in W^{\Lam}\circ \Lam$,
by
\begin{align*}
\ch V=\sum_{\mu\in W^{\Lam}\circ \Lam}
[V:M(\mu)]\ch M(\mu).
\end{align*}
Recall
$$\ch M(\lam)=\frac{e^{\lam}}{
\prod_{\alpha\in \iroots_-}(1-e^{\alpha})^{\dim
\g_{\alpha}}\prod_{\alpha\in \nrroots}(1-e^{\alpha})
}.$$
\begin{Pro}\label{Pro:estimate-of-character-revised}
Let $\Lam\in \dual{\h}_{\kappa}$
be as in Lemma \ref{Lem:estimate-height}.
Then,
for any $V\in
\BGG_{\kappa}^{[\Lam]}$,
\begin{align*}
&\sum_{\lam\in \dual{\h}}\ch H^{\bullet}
(\B_{\lam}^{\bullet}(\n,V))
e^{\lam}\\
&
\leq
\prod_{\alpha\in \iroots_-}(1-e^{\alpha})^{\dim
\g_{\alpha}}
\ch (V/\m V)
\prod\limits_{\alpha
\in \rroots(\n_-)\backslash \rroots(\m)}(1+z\inv e^{\alpha})
\prod\limits_{\alpha
\in \rroots(\n^t_-)}(1+z e^{\alpha})\\
&=\frac{\prod\limits_{\alpha
\in \rroots(\n_-)\backslash \rroots(\m)}(1+z\inv
e^{\alpha})
\prod\limits_{\alpha
\in \rroots(\n^t_-)}(1+z e^{\alpha})
}
{\prod\limits_{\alpha
\in  \nrroots
\backslash \rroots(\m)
}(1-e^{\alpha})}
\sum_{\mu\in W^{\Lam}\circ
\Lam}[V:M(\mu)]e^{\mu}
\end{align*}
where inequity $\leq $ means
that each coefficient of $z^ie^{\lam}$
of the left-hand-side is smaller than
or equal to
that of the right-hand-side.
\end{Pro}
\begin{proof}
Consider the
(obvious semi-infinite analogue of) Hochschild-Serre
spectral sequence for the subalgebra $\m\subset \n
$.
It is easy to check that
the corresponding filtration
is bounded upper
on each $C(\n,V)^{\lam}$,
$\lam\in \dual{\h}$.
By
definition,
\begin{align*}
E_1^{p,q}
=H_{-q}\left(\m,V\*
\Lam^{\semiinf +p}\left(\n/\m\right)
\right).
\end{align*}
Here,
$\Lam^{\semiinf +p}\left(\n/\m\right)
=\sum\limits_{i-j=p}
\Lam^i \n^t_-
\* \Lam^j \left(\n_-/\m\right)$
and $\m$ acts on $\Lam^i \n^t_-$
via the identification
$\n^t_-=\dual{\n}_+=\dual{(\n/\n_-)}$.
Clearly,
we have
\begin{align}\label{eq:enequlity-}
\ch \coho{\bullet}(\n,V)
\leq \sum_{p,q}z^{p+q} \sum_{\lam}
e^{\lam}\dim H_{-q}\left(\m,V\*
\Lam^{\semiinf +p}\left(\n/\m\right)
\right)^{\lam}.
\end{align}
Since $V$ is a free
$\m$-module, so is $V\*
\Lam^{\semiinf +p}\left(\n/\m\right)$,
$p\in \Z$.
Thus,
\begin{align}\label{eq:E_1:setimate-1}
E_1^{p,q}=\begin{cases}
\left(V\*
\Lam^{\semiinf +p}\left(\n/\m\right)\right)
/\m&
(q=0)\\
\{0\}&(q\ne 0).
\end{cases}
\end{align}
By \eqref{eq:enequlity-}
and \eqref{eq:E_1:setimate-1},
we get
\begin{align*}
\ch \coho{\bullet}(\n,V)
\leq
\ch (V/\m V)\cdot  \ch\Lam^{\semiinf
+\bullet}\left(\n/\m\right)\cdot
\prod_{\alpha
\in
\rroots(\m)}(1-e^{\alpha}).
\end{align*}
Here,
we have set
$ \ch\Lam^{\semiinf
+\bullet}\left(\n/\m\right)
=\sum_iz^i \sum_{\lam}e^{\lam}\dim
\left(\Lam^{\semiinf
+i}\left(\n/\m\right)\right)^{\lam}$.
It is easy to see that
\begin{align*}
\ch\Lam^{\semiinf
+\bullet}\left(\n/\m\right)
=\prod\limits_{\alpha
\in \rroots(\n_-)\backslash \rroots(\m)}(1+z\inv e^{\alpha})
\prod\limits_{\alpha
\in \rroots(\n^t_-)}(1+z e^{\alpha}).
\end{align*}
Therefore,
\eqref{eq:ch-B-and-Ch-C}
and \eqref{eq:ch-of-m-coinvariants}
%the facts that
%$\rroots(\t_{\srho\che}(\n_-))=
%\rroots(\n_-)\backslash \rroots(\m)$
%and
%$\t_{\srho\che}(\nrroots)
%\cap \nrroots=\rroots(\t_{\srho\che}(\n_-))\sqcup
%\rroots(\n_-^t)$
prove the Proposition.
\end{proof}

\begin{proof}[Proof of Proposition
\ref{Pro:pre-crutial-estimate}]
Suppose
$H^{i}(B_{\lam}(\n,V))\ne \{0\}$
for some $\lam\in \dual{\h}$.
Since
$\rroots(\m)=\w_{\n}(\prroots)
\cap \nrroots$,
we have
$\nrroots\backslash \rroots(\m)
=\nrroots\cap \w_{\n}(\nrroots)
=-\prroots\cap \w_{\n}(\prroots)$.
Therefore,
by Proposition
\ref{Pro:estimate-of-character-revised},
$\lam$ has the form as
\begin{align}\label{equation;form-of-weight}
\lam=\mu-\sum\limits_{\alpha\in
\prroots\cap \w_{\n}(\prroots)
}m_{\alpha}\alpha,
\quad \text{with }\sum m_{\alpha}\geq |i|,
\end{align}
with $\mu\in
W^{\Lam}\circ \Lam$
such that
$[V:M(\mu)]\ne 0$.
We claim that
\eqref{equation;form-of-weight}
implies
\begin{align}\label{eq:revised1-estimate}
\Whw{\n}{\lam}
\leq \Whw{\n}{\mu}-|i|.
\end{align}
Indeed,
for the $\n=\Ln{-}$ case
\eqref{eq:revised1-estimate}
easily  follows
from
\eqref{eq:diffrence-of-weights-}
and the fact that
$\prroots\cap w_0(\prroots)\cap
\sproots=\emptyset$.
To see \eqref{eq:revised1-estimate}
for the $\n=\Ln{+}$ case,
notice that
$t_{-\srho\che}(\prroots)\cap \sproots=\emptyset$,
and,
therefore,
$\bra t_{-\srho\che}(\alpha),\Dg\ket\geq 1$
for any
$\alpha\in \prroots\cap t_{\srho\che}(\prroots)$.
Then,
\eqref{eq:revised1-estimate}
follows from
\eqref{eq:diffrence-of-weights}.

By Proposition \ref{Pro:pre-crutial-estimate},
we have shown
that
$\DW_{\n}$-eigenvalues of
$\Wcoho^{i}\left(\n,V\right)$
is contained in the set
$
\bigcup\limits_\ud{\mu\in
W^{\Lam}\circ \Lam
}{
[V:M(\mu)]\ne 0} \Whw{\n}{\mu}-\Z_{\geq
|i|}$.
But
$[V:M(\mu)]\ne 0$
implies there exists
$\mu'\in W^{\Lam}\circ \Lam$
such
that
$[V:L(\mu')]\ne 0$
and $\mu'-\mu\in Q_+^{\Lam
}$.
Thus,
Proposition follows from
Lemma \ref{Lem:estimate-height}.
\end{proof}

\section{Vanishing of cohomology}
\label{section:Main-theorem}
\subsection{Vanishing of  cohomology
associated to projective modules
and injective modules}
For a given
$\Lam\in \dual{\h}_{\kappa}$,
let $\BGG_{\kappa}^{[\leq \Lam]}$
be the
full subcategory of
$\BGG_{\kappa}^{[\Lam]}$ consisting of module $V$ such
$V^{\lam}=\{0\}$
unless $\lam\in \Lam-Q_+$.
Then,
every
finitely generated object
of $\BGG_{\kappa}^{[\leq \Lam]}$
is an image of some projective object
of $\BGG_{\kappa}^{[\leq \Lam]}$
by \cite[2.10]{Moody}.
Let ${}^{\Delta}\BGG_{\kappa}^{[\leq \Lam]}$
be the full subcategory of
$\BGG_{\kappa}^{[\leq \Lam]}$
consisting of modules $V$
that admits a
Verma flag,
i.e,
a finite filtration
\begin{align*}
V=V_0\supset V_1\supset \dots \supset V_k=\{0\}
\end{align*}
such that each successive subquotient $V_i/V_{i+1}$
is isomorphic to some
Verma module.
It is known that
an object $V$ in $\BGG_{\kappa}^{[\leq \Lam]}$
belongs to
${}^{\Delta}\BGG_{\kappa}^{[\leq \Lam]}$
if and only if $\Ext_{\BGG_{\kappa}}^1(V,M(\lam)^*)=\{0\}$
for all $M(\lam)^*\in \BGG_{\kappa}^{[\leq \Lam]}$.
In particular,
projective objects in $\BGG_{\kappa}^{[\leq \Lam]}$ are
objects in ${}^{\Delta}\BGG_{\kappa}^{[\leq \Lam]}$.
\begin{Th}\label{Th:vanishing_of_projectives}$ $
For a given $\Lam\in \dual{\h}_{\kappa}$,
$\Wcoho^{i}(\n,V)=\{0\}$
$(i\ne 0)$ for all $V\in {}^{\Delta}\BGG_{\kappa}^{[\leq
\Lam]}$. In particular,
$\Wcoho^{i}(\n,P)=\{0\}$ $(i\ne 0)$
for all projective objects in $\BGG_{\kappa}^{[\leq
\Lam]}$.
\end{Th}
\begin{proof}
%Let $k(V)$ be
%the length $k$
We prove
by induction on the length $k(V)$
of the
Verma flag of $V$.
We have already proved the $k(V)=1$ case in
Theorem \ref{Th:vanishing.for.Verma}.
Let $k(V)\geq 2$.
Then,
there exits
an exact sequence
$0\rightarrow V_1\rightarrow V\rightarrow
M(\mu)\rightarrow 0$
($\mu\in \dual{\h}_{\kappa}$)
%\quad (V_1\in \BGG_{\kappa}^{\Delta},\lam\in \dual{\h}_{\kappa}),
in ${}^{\Delta}\BGG_{\kappa}^{[\leq \Lam]}$.
Thus,
the
corresponding long exact sequence
and the induction hypothesis
prove the proposition.
\end{proof}

Similarly,
let ${}^{\nabla}\BGG_{\kappa}^{[\leq \Lam]}$
be the full subcategory of
$\BGG_{\kappa}^{[\leq \Lam]}$
consisting of modules $V$
such that
$V^*\in{}^{\Delta}\BGG_{\kappa}^{[\leq \Lam]}$.
\begin{Th}$ $\label{Th:vanishing-injective}
Let $\Lam\in \dual{\h}_{\kappa}$ such that
\begin{align*}
&\begin{cases}
\text{$\bra \Lam+\rho,\alpha\che\ket\not\in \Z$
for all
$\alpha\in \prroots\cap t_{\srho\che}(\nrroots)$}
&(\text{if $\n=\Ln{+})$},
\\
\text{$\bra \Lam+\rho,\alpha\che\ket\not\in \Z$
for all
$\alpha\in \sproots$}
&(\text{if $\n=\Ln{-})$}.
\end{cases}
\end{align*}
Then,
$\Wcoho^{i}(\n,V)=\{0\}$
$(i\ne 0)$ for all $V\in
{}^{\nabla}\BGG_{\kappa}^{[\leq \Lam]}$.
In particular,
$\Wcoho^{i}(\n,I)=\{0\}$ $(i\ne 0)$
for all injective objects in
$\BGG_{\kappa}^{[\leq \Lam]}$.
\end{Th}
\begin{proof}
The assumption on $\Lam$
implies
$\Wcoho^{i}(\n,M(\lam)^*)=\{0\}$
$(i\ne0)$
for $\lam\in W^{\Lam}\circ \Lam$
by Theorem \ref{Th:vanishin-dual}.
Thus,
the theorem can be proved similarly
as Theorem \ref{Th:vanishing_of_projectives}.
\end{proof}
\subsection{Main theorem}
\begin{Th}\label{Theorem:Main-Theorem}
Let $\kappa\in \C\backslash \{0\}$
and $\Lam\in \dual{\h}_{\kappa}$.
\begin{enumerate}
\item
Suppose that
$\bra \Lam+\rho,\alpha\che\ket\not\in \Z$
for all
$\alpha\in \prroots\cap t_{\srho\che}(\nrroots)$.
Then,
\begin{align*}
\Wcoho^{i}(\Ln{+},V)=\{0\}\quad (i\ne0)
\end{align*}
for all
$V\in \BGG_{\kappa}^{[\Lam]}$.
\item
Suppose that
$\bra \Lam+\rho,\alpha\che\ket\not\in \Z$
for all
$\alpha\in \sproots$.
Then,
\begin{align*}
\Wcoho^{i}(\Ln{-},V)=\{0\}\quad (i\ne 0)
\end{align*}
for all
$V\in \BGG_{\kappa}^{[\Lam]}$.
\end{enumerate}
\end{Th}
\begin{proof}[Proof of Theorem
\ref{Theorem:Main-Theorem} when $\kappa
\in \C\backslash \Q_{\geq 0}$]

We may assume that
$V\in \BGG_{\kappa}^{[\leq\Lam]}$.
Since $\kappa
\in \C\backslash \Q_{\geq 0}$,
the cohomological dimension of $V\in \BGG_{\kappa}^{
[\leq \Lam]}$ is finite,
that is,
there exists a projective resolution
$$0\rightarrow P_n\rightarrow P_{n-1}\rightarrow \dots
\rightarrow P_0\rightarrow V\rightarrow 0$$
of $V$ in
$\BGG_{\kappa}^{[\leq \Lam]}$.
Let $N_k=\im \partial_k$.
Then,
$0\rightarrow N_{k+1}\rightarrow P_{k}\rightarrow N_{k}\rightarrow 0$.
Thus,
by the long exact sequence of semi-infinite cohomology,
we get
$H^i(N_{k})\cong  H^{i+1}(N_{k+1})$ for $i>0$
by Theorem \ref{Th:vanishing_of_projectives}.
This implies $H^i(V)=H^{i+n}(P_n)=\{0\}$
for all $i>0$.
The proof of $H^i(V)=\{0\}$ for $i<0$ is similar.
\end{proof}
When $\kappa\in \Q_{>0}$,
some modification of the proof is needed:
\begin{Pro}\label{Pro:estimate-of-shifting}
Let $\kappa\in \Q_{>0}$
and
let $\Lam\in \dual{\h}_{\kappa}$
be as in Theorem \ref{Th:vanishing-injective}.
Let $V\in
\BGG_{\kappa}^{[\Lam]}$.
\begin{enumerate}
\item
For a given $a\in \C$,
there exist an object $N$ in $\BGG_{\kappa}^{[ \Lam]}$
such that
$\Wcoho^{i}(\n,V)_a=
\Wcoho^{i+1}(\n,N)_a$
for all $i>0$.
\item
For a given $a\in \C$,
there exist an object $N'$ in $ \BGG_{\kappa}^{[ \Lam]}$
such that  $\Wcoho^{i}(\n,V)_a=
\Wcoho^{i-1}(\n,N')_a$
for all $i<0$.
\end{enumerate}
\end{Pro}
\begin{proof}
(1) By Proposition
\ref{Pro:crutial-estimate},
there exists
finitely generated
submodule $V'$ of $V$
such that
\begin{align}\label{eq:in-proof-inter-iso}
\Wcoho^{i}(\n,V)_a
\cong  \Wcoho^{i}(\n,V')_a.
\end{align}
for the given $a$.
Since
$V'$ is finitely generated,
there exists
some projective
object $P$
of $\BGG_{\kappa}^{[
\Lam]}$
and
an exact sequence
$0\rightarrow N\rightarrow P\rightarrow
V'\rightarrow 0$
in $\BGG_{\kappa}^{[
\Lam]}$.
Therefore, we get
$\Wcoho^{i}(\n,V')\cong
\Wcoho^{i+1}(\n,N)$ for
all $i>0$ by
Theorem \ref{Th:vanishing_of_projectives}.
By \eqref{eq:in-proof-inter-iso},
this implies
$\Wcoho^{i}(\n,V)_a
\cong  \Wcoho^{i+1}(\n,N)_a$
for all $i>0$.
(2) can be similarly proved
by using Theorem \ref{Th:vanishing-injective}.
\end{proof}
\begin{proof}[Proof of Theorem
\ref{Theorem:Main-Theorem} when $\kappa
\in \Q_{\geq 0}$]
We may assume $\Lam\in \fd_{\kappa.+}$.
It is sufficient
to show
that
$\Wcoho^{i}(\n,V)_a=\{0\}$
($i\ne 0$)
for all $V\in \BGG_{\kappa}^{[\Lam]}$
and $a\in \C$.

Fix $a\in \C$.
By applying Proposition
\ref{Pro:estimate-of-shifting} (1) repeatedly,
it follows that,
for any $r>0$,
there exists an object
$N_r$
of $\BGG_{\kappa}^{[
\Lam]}$
such that
\begin{align*}
\Wcoho^{i}(\n,V)_a\cong
\Wcoho^{i+r}(\n,N_{r})_a
\quad (i>0).
\end{align*}
This forces
$\Wcoho^{i}(\n,V)_a=\{0\}$
for $i>0$
by
Corollary \ref{CO:estimate-of-weights-kappa-Q-+}.
The proof for  $i<0$ is similar.
\end{proof}
The
following is straightforward
from  Theorem
\ref{Theorem:Main-Theorem}
and  Remark
\ref{Rem:character-of-VErma}.
\begin{Co}\label{Co:exact-fun}
Let $\Lam$ be as in Theorem
\ref{Theorem:Main-Theorem}.
Then,
the correspondence
$V\rightsquigarrow \Wcoho^{0}(\n,V)$
defines an exact functor
from $\BGG_{\kappa}^{[\Lam]}$
to the category of $\W_{\kappa}(\g)$-modules.
In particular,
%$\dim \Wcoho^{0}(\n,V)_a <\infty$
%for all $a\in \C$
%and
\begin{align*}
\ch \Wcoho^{0}(\n,V)=
\sum_{\mu\in W^{\Lam}\circ \Lam}
[V:M(\mu)]\frac{q^{\Whw{\n}{\mu}}}{
\prod_{i\geq 1}(1-q^i)^{\rank \sg}}
\end{align*}
for $V\in \BGG_{\kappa}^{[\Lam]}$.
\end{Co}

\begin{Rem}\label{Rem:final}$ $

\begin{enumerate}
\item Let $\kappa\in \C\backslash \Q$.
Then,
any
$\Lam\in \dual{\h}_{\kappa}$
such that $\bar{\Lam}\in \sP$
satisfies the condition
of Theorem \ref{Theorem:Main-Theorem} (1)
\item\label{Rem:vanish-FKW}
It was proved in \cite{FKW} that
$ \Wcoho^{\bullet}(\Ln{-},L(\lam))\equiv 0$
if $\bra \lam,\alpha\che\ket\in \Z_{\geq 0}$
for some $\alpha\in \sPi$.
%(cf. \cite[Proposition 6.4]{Backelin}).
\item Suppose that
$\bra \Lam+\rho,\alpha\che\ket\not\in \Z$
for all
$\alpha\in \sproots$.
Then,
by Corollary \ref{Co:exact-fun} and Remark
\ref{Rem:exact-twist}, it follows
that
\begin{align*}
\ch \Wcoho^{0}(\Ln{-},L(\lam))=
\ch \Wcoho^{0}(\Ln{-},L(w\circ \lam))
\end{align*}
for $w\in \sW$ and $\lam\in W^{\Lam}\circ \Lam$.
\item
A  principal admissible weight $\Lam$
(\cite{KW2}) is called {non-degenerate} if
$\bra \Lam
+\rho,\alpha\che \ket\not\in \Z$
for all $\alpha\in \sproots$.
By Theorem \ref{Theorem:Main-Theorem} (2),
it follows that
$\Wcoho^{i}(\Ln{-},L(\Lam))=\{0\} $
($i\ne 0$)
if $\Lam$ is a non-degenerate
principal admissible weight.
This was conjectured
by Frenkel-Kac-Wakimoto
(\cite[Conjecture 3.4$_-$]{FKW}).
%is proved for the $\Ln{-}$ case.
\item\label{rem-on-minimal}
Let $\kappa=p/q$,
$p\in \Z_{ \geq h\che}$,
$q\in \Z_{\geq h}$,
$(p,q)=1$,
$(q,r\che)=1$,
where $h$ is the Coxeter number of $\sg$
and $r\che=\max\{ 2/(\alpha,\alpha);\alpha\in \sPi\}$.
Then,
$\kappa-h\che$ is
a principal admissible number
(\cite{KW2}).
Set
\begin{align*}
\Lam_{\lam,\mu}=\lam-\kappa \mu+(\kappa-h\che)\Lam_0
\quad
\text{($\lam\in \bar{P}^{p-h\che}_+$,
$\mu \in \bar{P}^{\vee q-h}_+$)},
\end{align*}
where
$\bar{P}^{p-h\che}_+=\{\lam\in  \bar{P}\mid
0\leq \bra \lam,\alpha\che\ket \leq p-h\che~
(\forall \alpha\in \sproots)\}$
and
$\bar{P}^{\vee q-h}_+=\{\mu\in \bar{P}^{\vee}\mid
0\leq \bra \alpha,\mu\ket \leq q-h~
(\forall \alpha\in \sproots)\}$.
Let
\begin{align*}
\Dot{\Pr}=\{\Lam_{\lam,\mu};
(\lam,\mu)\in\bar{P}^{p-h\che}_+
\times \bar{P}^{\vee q-h}_+
\}\subset \dual{\h}_{\kappa}.
\end{align*}
Then,
$\Dot{\Pr}$
is a subset of
the set of
principal admissible weights
of $\g$ of level $\kappa-h\che$.
Note that
$(\kappa-h\che)\Lam_0=\Lam_{0,0}\in \Dot{\Pr}$.

It is easy to see that
any element of
$\Dot{\Pr}$
satisfies the condition
of Theorem \ref{Theorem:Main-Theorem} (1).
Thus,
\begin{align}\label{eq_in_Remark:vanishing}
\Wcoho^i(\Ln{+},L(\Lam))=\{0\}
\quad (i\ne 0)
\quad \text{for
$\Lam \in \Dot{\Pr}$}
\end{align}
by Theorem \ref{Theorem:Main-Theorem} (1).
This proves
the
conjecture
of Frenkel-Kac-Wakimoto
\cite[Conjecture 3.4$_+$]{FKW}
partially.
Note that
\eqref{eq_in_Remark:vanishing}
in particular implies
\begin{align*}
\Wcoho^i(\Ln{+},L((\kappa-h\che)\Lam_0))=\{0\}
\quad (i\ne 0).
\end{align*}
%Hence,
%$\Wcoho^0(\Ln{+},L((\kappa-h\che)\Lam_0))$
%is a vertex operator algebra
% by the result of \cite{FM}.
It is expected  that
$\Wcoho^0(\Ln{+},L((\kappa-h\che)\Lam_0))$
is a rational VOA
and
that
the modules
$\{\Wcoho^0(\Ln{+},L(\Lam));
\Lam\in \Dot{\Pr}\}$
exhaust
the simple objects
of the vertex
operator
algebra
$\Wcoho^0(\Ln{+},L((\kappa-h\che)\Lam_0))$, see
\cite{FKW}.
\end{enumerate}\end{Rem}

\end{document}